\documentclass[reqno,10pt]{amsart}
\usepackage{xcolor}

\usepackage{amssymb}
\usepackage{comment}

\usepackage[all,cmtip]{xy}
\newtheorem{MainTheorem}{Theorem}
\newtheorem{Proposition}{Proposition}[section]
\newtheorem{Definition}[Proposition]{Definition}
\newtheorem{Lemma}[Proposition]{Lemma}
\newtheorem{Theorem}[Proposition]{Theorem}
\newtheorem*{Question}{Question}
\newtheorem*{Theorem*}{Theorem}
\newtheorem{Corollary}[Proposition]{Corollary}
\newtheorem{Remark}[Proposition]{Remark}

\DeclareMathOperator{\Val}{Val}

\DeclareMathOperator{\Gr}{Gr}
\DeclareMathOperator{\AGr}{AGr}
\DeclareMathOperator{\Cr}{Cr}
\DeclareMathOperator{\Res}{Res}

\DeclareMathOperator{\Id}{Id}
\DeclareMathOperator{\Jac}{Jac}
\renewcommand{\Re}{\mathrm {Re}}
\renewcommand{\Im}{\mathrm {Im}}
\DeclareMathOperator{\vol}{vol}

\DeclareMathOperator{\Dens}{Dens}
\DeclareMathOperator{\Span}{Span}

\DeclareMathOperator{\Stab}{Stab}
\DeclareMathOperator{\Sym}{Sym}
\DeclareMathOperator{\Ker}{Ker}
\DeclareMathOperator{\tr}{tr}

\DeclareMathOperator{\Kl}{Kl}

\DeclareMathOperator{\WF}{WF}

\DeclareMathOperator{\supp}{supp}

\DeclareMathOperator{\GL}{GL}
\DeclareMathOperator{\SL}{SL}
\DeclareMathOperator{\sign}{sign}
\DeclareMathOperator{\ev}{ev}

\newcommand{\R}{\mathbb{R}}
\newcommand{\K}{\mathcal{K}}
\newcommand{\C}{\mathbb{C}}
\newcommand{\OO}{\mathrm{O}}
\newcommand{\SO}{\mathrm{SO}}

\makeatletter
\newcommand{\appendixauthor}[1]{%
	{\parindent0pt\vspace*{-25pt}%
		\linespread{1.1}\large\scshape#1%
		\par\nobreak\vspace*{35pt}}
	\@afterheading%
}
\makeatother
\title{Crofton formulas and indefinite signature}
\author{Dmitry Faifman}

\email{dfaifman@math.utoronto.ca}

\address{Department of Mathematics, University of Toronto, Bahen Centre, 40 St. George St., Toronto ON M5S 2E4, Canada}

\thanks{This research was partially supported by an NSERC Discovery grant}

\date{}

\begin{document}

\begin{abstract}	
	We study the $\OO(p,q)$-invariant valuations classified by A. Bernig and the author. Our main result is that every such valuation is given by an $\OO(p,q)$-invariant Crofton formula. This is achieved by first obtaining a handful of explicit formulas for a few sufficiently general signatures and degrees of homogeneity, notably in the $(p-1)$ homogeneous case of $O(p,p)$, yielding a Crofton formula for the centro-affine surface area when $p\not\equiv 3\mod 4$. We then exploit the functorial properties of Crofton formulas to pass to the general case. We also identify the invariant formulas explicitly for all $\OO(p,2)$-invariant valuations. The proof relies on the exact computation of some integrals of independent interest. Those are related to Selberg's integral and to the Beta function of a matrix argument, except that the positive-definite matrices are replaced with matrices of all signatures.
	We also analyze the distinguished invariant Crofton distribution supported on the minimal orbit, and show that, somewhat surprisingly, it sometimes defines the trivial valuation, thus producing a distribution in the kernel of the cosine transform of particularly small support.
	In the heart of the paper lies the description by Muro of the $|\det X|^s$ family of distributions on the space of symmetric matrices, which we use to construct a family of $\OO(p,q)$-invariant Crofton distributions. We conjecture there are no others, which we then prove for $\OO(p,2)$ with $p$ even. 
	The functorial properties of Crofton distributions, which serve an important tool in our investigation, are studied by T. Wannerer and the author in the Appendix.
	
\end{abstract}
\maketitle

\section{Introduction}

\subsection{Overview.}

This paper deals with Crofton formulas, which lie within the domain of integral geometry, and have numerous applications in convex and stochastic geometry. For an exposition of those topics, see \cite{klain_rota} and \cite{schneider_stochastic}. Crofton formulas can be considered as the simplest instance of kinematic formulas, a central theme in integral geometry that has recently seen intensive development in the framework of convex valuation theory.

A presentation of the classical theory of valuations, starting with Buffon's needle problem and Dehn's solution of Hilbert's 3rd problem, and spanning contributions by Minkowski, Alexandrov, Blaschke, Santal\'{o}, Hadwiger, McMullen, Schneider and others, can be found in \cite{schneider_book14}, \cite{klain_rota} and the references therein. More recently - see \cite{alesker_survey}, \cite{bernig_survey}, \cite{fu_survey} for recent surveys - a rapid progress in integral geometry followed Alesker's solution of McMullen's conjecture \cite{alesker_mcmullen}, allowing also to relate the theory to past contributions by other geometers, such as Weyl's tube formula, Chern's kinematic formulas and more, as well to the other, Gelfand-style branch of integral geometry, studying Radon transforms and their generalizations, see \cite{alesker_intgeo}. Valuation theory has also been recently used to obtain new types of Brunn-Minkowski \cite{log_concave} and Alexandrov-Fenchel \cite{abardia_wannerer} inequalities.

The classical Crofton formula in its simplest form computes the length of a rectifiable compact curve $\gamma\subset \R^2$: 
\[\text{Length}(\gamma)=\frac12\int_{\overline{\R P}^1}\chi(\overline L\cap \gamma)d\overline L \]
where $\overline{\R P}^1$ is the 2-dimensional manifold of affine lines in $\R^2$ equipped with the appropriately normalized rigid motion invariant measure $d\overline L$, and $\chi$ is the Euler characteristic. This formula can be used to produce a simple proof of the isoperimetric inequality in the Euclidean plane.

More generally, a Crofton formula can be written for all Euclidean intrinsic volumes (also known as quermassintegrals). Those include the surface area $\mu_{n-1}$ and the mean width $\mu_1$. For $1\leq k\leq n-1$ and $M\subset \R^n$ a sufficiently nice compact subset (e.g. a convex set, or a submanifold with corners) one has
\[\mu_k(M)=c_{n,k}\int_{{\AGr}_{n-k}(\R^n)}\chi(M\cap \overline E)d\overline E  \]
where the integral is taken over the Grassmannian of affine $(n-k)$-subspaces, with respect to the rigid-motion invariant measure, and $c_{n,k}$ is an explicit constant. 

Those formulas fit neatly into the general framework of convex valuation theory. Put simply, a convex valuation on $V=\R^n$ is a finitely additive measure on compact convex sets (henceforth: the convex bodies). Of particular interest are the elements of $\Val(V)$, the translation-invariant valuations that are continuous with respect to the Hausdorff metric on convex bodies. All valuations in this paper are translation-invariant, and we will usually omit the term.

We have the Crofton map
\[\Cr:\mathcal M(\AGr_{n-k}(V))^{tr} \to \Val(V) \]
from the space of translation-invariant measures on the affine Grassmannian to the space of translation-invariant continuous valuations, given by
\[\Cr(\mu)(K)=\int_{{\AGr}_{n-k}(V)}\chi(K\cap \overline E)d\mu(\overline E) \]
One immediately sees that the image of $\Cr$ falls inside the space of even $k$-homogeneous valuations $\Val_k^+(V)=\{\phi\in \Val(V):\phi(\lambda K)=|\lambda|^k\phi(K), \forall\lambda\in\R,K\in \K(V)\}$. We say that the corresponding valuations are given by a Crofton formula, and refer to the elements of $\mathcal M(\AGr_{n-k}(V))^{tr} $ as $k$-homogeneous Crofton measures.

In \cite{alesker_bernstein}, Alesker and Bernstein have shown that the image of  $\Cr$ is dense in $\Val_k^+(V)$, equipped with the topology of uniform convergence on compact sets. In a later paper \cite{alesker_faifman}, Alesker and the author have shown that in fact, any valuation $\phi\in \Val_k^+(V)$ can be represented as $\Cr(\mu)$, where $\mu\in \mathcal M^{-\infty}(\AGr_{n-k}(V))^{tr}$ is some translation-invariant distribution (generalized measure). Moreover, if one replaces the even continuous valuations with the somewhat larger class of even generalized valuations $\Val^{+,-\infty}(V)$ (in which the former is a dense subspace), one can extend $\Cr: \mathcal M^{-\infty}(\AGr_{n-k}(V))^{tr} \to \Val_k^{+,-\infty}(V)$
as a surjection.

Given a group $G\subset \GL(V)$, write $X^G$ for the $G$-invariants in a $G$-module $X$. A natural problem is to describe the space $\Val(V)^G$ of $G$-invariant valuations on $V$. For $G=\OO(n)$, Hadwiger's theorem \cite{hadwiger} states that the invariant valuations are precisely the intrinsic volumes $\mu_k$. Generalizing Hadwiger's results, Alesker has shown that for a compact Lie group acting transitively on the space of lines $\mathbb P(V)$, $\dim \Val(V)^G<\infty$, and moreover any $\phi\in \Val(V)^G$ is given by $\phi=\Cr(\mu)$ with a smooth Crofton measure $\mu$. 

In recent years, explicit Hadwiger-type results were obtained for the various compact groups $G$ with this property. Most notably, the valuation theory of the complex unitary group  $G=\mathrm U(n)$ was understood completely - an explicit basis of invariant valuations was given by Alesker \cite{alesker03_un}, and subsequently the full array of kinematic formulas was determined by Fu \cite{fu06} and Bernig-Fu \cite{bernig_fu_hig}. In particular, a full array of Crofton formulas is available for the $\mathrm U(n)$-invariant valuations. Other groups were considered in \cite{bernig_g2}, \cite{bernig_su}, \cite{bernig_solanes}, \cite{bernig_voide}.
\\\\
For a general non-compact Lie group $G$, less is known. The non-trivial invariant valuations typically exhibit one type of discontinuity or another. In this note, we are concerned with $G$-invariant generalized valuations, which are one natural extension of continuous valuations, and can be thought of as valuations on smooth bodies. A different approach is taken by Ludwig and Reitzner \cite{ludwig_reitzner99, ludwig_reitzner10}, who classify the upper semi-continuous valuations invariant under the action of $G=\SL(n)$, with or without translation-invariance; in the former case, the only non-obvious invariant is the affine surface area. 

\subsection{Main results.}
Assume $-\Id\in G$, so all invariant valuations must be even. As we mentioned, it is a corollary of Alesker's irreducibility theorem that every even (generalized) valuation is given by a Crofton formula.\\ 
\begin{Question}
	Is every $G$-invariant valuation given by a $G$-invariant Crofton formula?
\end{Question}

The answer is obviously positive for compact $G$, through averaging over $G$. The main goal of this paper is to give a positive answer for $G=\OO(p,q)$, the symmetry group of a non-degenerate quadratic form of indefinite signature, arguably the simplest family of non-compact Lie groups. We also obtain some partial results towards an explicit description of those formulas.

Despite the superficial similarity to the Euclidean case, the case of indefinite signature is much harder to study. 
One obvious reason is that there is no natural compact body attached to the group. Another reason is that smooth measures get replaced with distributions, which are often defined indirectly through a meromorphic extension of an integral converging elsewhere. Other difficulties arise when the maximal compact subgroup becomes too small in some sense, which happens when $\min(p,q)\geq 2$.
In the Euclidean case, one can essentially restrict attention to the combinatorics of dissections of polytopes, as was done by Hadwiger. The indefinite signature however is inherently adapted for smooth convex bodies. The integral geometry of $\OO(p,q)$ brings together such diverse subjects as the representation theory of symmetric spaces, microlocal analysis and matrix integrals.

The Lorentz group $\OO(n-1,1)$ was considered by S. Alesker and the author in \cite{alesker_faifman}, and the general signature was studied by A. Bernig and the author in \cite{bernig_faifman_opq}. There, the dimensions of the spaces of invariant valuations were computed, and a simple description was given in terms of their Klain sections. 

In \cite{bernig_faifman_opq}, a complete set of Crofton formulas was obtained for $\R^{2,2}$, while in \cite{alesker_faifman} a complete set of Crofton distributions was constructed for the Lorentz signature $(n-1,1)$. In both of those cases, the Crofton distributions could be chosen to be invariant. In this paper we consider the general case.

We first establish that the space of $\OO(p,q)$-invariant Crofton distributions has finite dimension. An upper bound on the dimension appears in Proposition \ref{prop:generalCroftonBound}.

Our first main result uses Muro's description of the meromorphic extensions of the distributions $|\det X|^s$ on spaces of symmetric matrices, to construct a family of $\OO(p,q)$-invariant Crofton distributions. 
A more precise statement, including fairly explicit descriptions of those Crofton distributions, appears in Theorem \ref{thm:1}.
\begin{MainTheorem}\label{thm:main_trivial}
There are at least as many linearly independent $\OO(p,q)$-invariant Crofton distributions on $\Gr_k(\R^{p,q})$ as there are open orbits under the action of $\OO(p,q)$, namely $\min(p,q,k,n-k)+1$. Equality is attained for $(p,q)=(2m,2)$.
\end{MainTheorem}
It is easy to check that equality holds for $\min(p,q)\leq 1$ as well. We conjecture that equality holds for general $(p,q)$. 
\\\\Write $\overline G$ for the group generated by $G\subset \GL(V)$ and all translations of $V$. The main result of the paper is the following.
\begin{MainTheorem}\label{thm:main_1}
	For all $n=p+q$ and $0\leq k\leq n$, $\Cr:\mathcal M^{-\infty}(\AGr_{n-k}(\R^{p,q}))^{\overline{\OO(p,q)}}\to  \Val^{-\infty}_k(\R^{p,q})^{\OO(p.q)}$ is surjective.
\end{MainTheorem}

The proof is surprisingly indirect. While we can write fairly explicit formulas for various $\OO(p,q)$-invariant Crofton distributions, it is not generally clear how to determine whether or not they in fact define non-trivial generalized valuations. While there is no a-priori reason to believe that in general there even exists one invariant valuation given by an invariant Crofton distribution, the opposite is certainly true, namely, there are invariant Crofton distributions defining the zero valuation. In fact, for every non-trivial degree of homogeneity, the space of invariant valuations is 2-dimensional, while the space of invariant Crofton distributions is generally larger by Theorem \ref{thm:main_trivial}, implying that for $2\leq k\leq n-2$, $p\geq q\geq 2$, the Crofton map on the invariant distributions has a non-trivial kernel. 

 Applying those Crofton formulas even to the simplest convex bodies appears to be computationally intractable in general, often resulting in a hypergeometric-like function of a matrix argument, with the positive matrices replaced by all matrices in the domain of integration. We are thus forced to evade computation. We first show that the existence of invariant Crofton formulas in all cases is implied by the existence of just one such formula in a sufficiently general infinite family of signatures and degrees of homogeneity (henceforth: universal family). This is done by exploiting the available symmetries of the $\OO(p,q)$-invariant valuations and their interconnectedness across all dimensions and signatures. 
 We then study certain invariant Crofton distributions in a few particularly amenable universal families of $(p,q,k)$.
 \begin{MainTheorem}\label{thm:main_universal_families} For $m\geq 1$ and $(p,q,k)$ either of: $(m, m, m-1)$, $(m, m, m+1)$, $(m+1, m, m)$, $(m+1, m, m+1)$, we find an explicit $\mu\in \mathcal M^{-\infty}(\AGr_{k}(\R^{p,q}))^{\overline{\OO(p,q)}}$ for which $\Cr(\mu)\neq 0$. \end{MainTheorem}
 We refer to Theorem \ref{thm:universal_families} for more details.
 
 The case of $(p,q)=(m,m)$ appears frequently in convex geometry, in particular it is directly related to the centro-affine surface area. Namely, for a convex set $K\subset V=\R^m$ containing the origin in its interior, consider $K^+:=\{(x,\xi)\in K\times K^o: x\cdot \xi =1\}\subset V\times V^*$. The space $V\times V^*$ has a natural quadratic form $Q(x,\xi)=x\cdot\xi$ of signature $(m,m)$. When $K$ is smooth and has positive gaussian curvature, $K^+$ is a smooth $(m-1)$-dimensional submanifold, and the restriction of $Q$ to $K^+$ is positive-definite. The centro-affine surface area is $\Omega_c(K)=\int_{K^+}\vol_{Q|_{K^+}}$. 
 
 \begin{MainTheorem}
 	
 	For $m\not \equiv 3\mod 4$ there is an explicit Crofton distribution $\mu_m\in \mathcal M^{-\infty}(\AGr_{m+1}(V\times V^*))^{\overline{\OO(m,m)}}$ such that
  	\[\Omega_c(K)=c_m\int_{\AGr_{m+1}(\R^{m,m})}\chi(K^+\cap\overline E)d\mu_m(\overline{E})\]
 	for some universal constant $c_m$.
 \end{MainTheorem}
For details see Theorem \ref{thm:centro_affine}.

The main technical tool for traveling between different $\R^{p,q}$ is that of restriction and projection of Crofton distributions, or more generally pull-back and push-forward of Crofton distributions by linear maps. This general construction is independent of the rest of the paper, and is coauthored by Thomas Wannerer.

While fairly straightforward in the smooth setting, it gets somewhat technically involved to extend to the case of distributions, which is essential in the $\OO(p,q)$-invariant setting, and we apply microlocal techniques such as wavefront set analysis.
\begin{MainTheorem}
	Let $f:U\to V$ be a linear map, $\dim V=n$ and $\dim U=n-d$, $d\in\mathbb Z$. There are partially defined natural maps $f_*$, $f^*$ between the corresponding spaces of Crofton distributions, such that the following diagrams commute whenever the maps are defined.
	\begin{displaymath}
	\xymatrix{
		\mathcal M^{-\infty}(\AGr_{n-k}(V))^{tr} \ar[r]^--{\Cr}\ar[d]^{f^*} &  \Val^{+,-\infty}_k(V) \ar[d]^{f^*} \\
		\mathcal M^{-\infty}(\AGr_{n-d-k}(U))^{tr} \ar[r]^--{\Cr} &  \Val^{+,-\infty}_k(U)}
	\end{displaymath}
	
	\begin{displaymath}
	\xymatrix{
		\mathcal M^{-\infty}(\AGr_{n-k}(U))^{tr}\otimes \Dens^*(U) \ar[r]^--{\Cr}\ar[d]^{f_*} &  \Val^{+,-\infty}_k(U) \otimes \Dens^*(U)\ar[d]^{f_*} \\
		\mathcal M^{-\infty}(\AGr_{n-k}(V))^{tr}\otimes \Dens^*(V) \ar[r]^--{\Cr} &  \Val^{+,-\infty}_{k-d}(V)\otimes \Dens^*(V)}
	\end{displaymath}
	The domains of definition are dense and contain the smooth Crofton measures.
\end{MainTheorem}
\noindent Here the vertical arrows on the right are the corresponding maps for valuations as defined by Alesker \cite{alesker_fourier}, and extended to generalized valuations by Bernig and the author in \cite{bernig_faifman_opq}. The precise statements can be found in Appendix \ref{sec:functorial}. 
\\\\
Under the action of $\OO(p,q)$ on $\Gr_k(\R^{p,q})$, there is a unique orbit of minimal dimension, which is also the unique closed orbit, denoted $X^k_c$. In about half the cases, there is a distinguished $\OO(p,q)$-invariant Crofton distribution $\mu_c$ supported on $X^k_c$. We show that it can define both a zero or a non-zero valuation in infinitely many cases. Namely, we establish the following.

\begin{MainTheorem}\label{thm:main_mu_c}
	\begin{enumerate}
		\item If $n\equiv \min(k, n-k, q)\mod 2$, there is a unique up to scale $\OO(p,q)$-invariant Crofton distribution, denoted $\mu_c$, which is supported on $X^k_c$. Otherwise, no such non-trivial Crofton distribution exists.
		\item For $q=2$, $p\geq 2$ even and arbitrary $2\leq k\leq p$, $\Cr(\mu_c)\neq 0$.
		\item For $m\geq 2$ and  $(p,q,k)=(2m, 2m-1, 2m-1)$, $\Cr(\mu_c)=0$.
		
	\end{enumerate}

\end{MainTheorem}
Part i) is just part iv) of Theorem \ref{thm:1}. Part ii) is contained in Theorem \ref{thm:q=2}. The last part is Theorem \ref{thm:mu_c}. We thus obtain a family of distributions lying in the kernel of the cosine transform with support of particularly large codimension.
\\\\
 Finally, in the case of signature $(p,2)$ we obtain the following.
\begin{MainTheorem}\label{thm:main_q=2}
	For all $p\geq 2$, there is an explicit basis of $\OO(p,2)$-invariant valuations given by $\OO(p,2)$-invariant Crofton distributions. Details appear in Theorem \ref{thm:q=2}.
\end{MainTheorem}

Again, instead of making explicit computations in a given space, we are forced to navigate between different dimensions to arrive at the result. 
\\\\
Theorem \ref{thm:main_universal_families} is proved by explicitly evaluating the corresponding Crofton formula on certain bodies. For $\R^{p,p}$ it is just the Euclidean ball, while for $\R^{p,p-1}$ it is the appropriately rescaled limit on a family of $\OO(p)\times \OO(p-1)$-symmetric ellipsoids degenerating to a $(p-1)$-dimensional ball. 
In all of those cases, the computation then boils down to the evaluation of certain linear combinations of matrix integrals. Namely, we study $D_n^\epsilon(s)=\displaystyle\sum_{a+b=n} \epsilon(b)\int_{S_{a,b}(1)}|\det X|^s dX$, where $S_{a,b}(1)=\{-I_n\leq X\leq I_n\}\cap S_{a,b}$ are the symmetric $n\times n$ real matrices of signature $(a,b)$; 
and $\epsilon(b)$ is any of the four coefficient functions: $\epsilon_{\textbf{abs}}(b)=1$, $\epsilon_{\textbf{sgn}}(b)=(-1)^b$, $\epsilon_{\textbf{cos}}(b)=\cos(\frac \pi 2 b)$,  $\epsilon_{\textbf{sin}}(b)=\sin(\frac \pi 2 b)$.

The positive-definite summand is an instance of the multivariate Beta function $B_n(a,b)=\int_{0\leq X\leq I_n}(\det X)^{a-\frac{n-1}{2}} \det(I-X)^{b-\frac{n-1}{2}}dX$ (which in turn is an instance of Selberg's integral). 
Our initial goal is to prove $D_n^\epsilon(s)$ does not vanish at a certain value $s=s_0$, on which it depends meromorphically ($s_0$ happens to fall outside the domain of convergence). Somewhat mysteriously, in all the cases we consider, the integral turns out to be non-vanishing at all values of $s$, leaving one wondering whether there is a deeper reason behind this phenomenon.

The following integral is just a multiple of $D_n^{\textbf{abs}}(s)=\int_{-I_n\leq X\leq I_n}|\det X|^sdX$. We denote $\Delta_n=\{\lambda\in\R^n: 1\geq \lambda_1\geq\dots\geq \lambda_n\geq -1\}$.
\begin{MainTheorem}\label{thm:Beta_integral} For $\Re(s)>-1$ one has
	\[ \displaystyle\int_{\Delta_n}  \prod_{i=1}^n |\lambda_i|^s \prod_{i<j}(\lambda_i-\lambda_j)d\lambda=  2^n \frac{ \displaystyle\prod_{\scriptsize\begin{array}{c}1\leq i< j\leq n\\i\equiv j\mod 2 \end{array} }(j-i) }   {\displaystyle\prod_{\scriptsize\begin{array}{c} 1\leq i\leq n\\i\equiv 1\mod 2 \end{array}} (s+i) \prod_{\scriptsize\begin{array}{c}1\leq i< j\leq n\\i\not\equiv j\mod 2 \end{array}} (2s+i+j)} \]
\end{MainTheorem}
Similar formulas hold for $D_{2n}^{\textbf{sgn}}(s)=\int_{-I_{2n}\leq X\leq I_{2n}}\sign(\det X)|\det X|^s dX$, as well as for $D_n^{\textbf{cos}}(s)$ and $D_n^{\textbf{sin}}(s)$. They can be found in section \ref{sec:selberg}.

All of those integrals loosely fall into, and are implied by, another peculiar family of integrals that we compute, akin to a family of integrals considered by Robbins \cite{robbins} and DiPippo-Howe \cite{dipippo_howe}.

\begin{MainTheorem}\label{thm:integral}
	
	For $e=(e_1,\dots, e_n)\in \mathbb N^n$ define the family of integrals \[f_n(e)=\int_{\Delta_n} \det (\lambda_i^{e_j-1})d\lambda\]
	Let $N_+(e)$ be the number of even entries of $e$, and $N_-(e)$ the number of odd entries. 
	\begin{enumerate}
		\item For a permutation $\sigma\in S_n$, $f_n(\sigma e)=(-1)^\sigma f_n(e)$.
		\item If $N_-(e)-N_+(e)\notin\{0,1\}$ then $f_n(e)=0$.
		\item Assume $N_-(e)-N_+(e)\in\{0,1\}$. The entries of $e$ can then be rearranged such that $e_i\equiv i \mod 2$. Then
		
		\[f_n(e)= (-1)^{n\choose 2} 2^n\frac {\displaystyle \prod_{\scriptsize\begin{array}{c}i< j\\i\equiv j\mod 2 \end{array}}(e_j-e_i) }{\displaystyle\prod_{i\equiv 1\mod 2} e_i \prod_{\scriptsize\begin{array}{c}i< j\\i\not\equiv j\mod 2 \end{array}} (e_i+e_j)}
		\]
		
	\end{enumerate}
\end{MainTheorem}
A similar formula is shown to hold when the $e_j$ are strict half-integers.%, see section \ref{sec:selberg} for details.

\subsection{Plan of the paper}
In section \ref{sec:preliminaries} we collect the necessary background from valuation theory and $O(p,q)$-geometry, and recall the classification of $\OO(p,q)$-valuations. In section \ref{sec:trivial_obstruction} we show by analyzing the symbols of the distribution that every orbit of $\Gr_k(\R^{p,q})$ contributes at most one dimension to the space of invariant Crofton distributions. It then remains to construct distributions corresponding to the invariant symbols on the different orbits, and analyze when can they be patched together to form a globally defined invariant Crofton distribution. In section \ref{sec:muro}, we show that a collection of such distributions can be obtained from a well-known family of homogeneous distributions on the space of symmetric matrices. In doing so we rely on results of Muro \cite{muro}, \cite{muro_prehomogeneous} who described those distributions very explicitly. This proves Theorems \ref{thm:main_trivial} and \ref{thm:main_mu_c} part i).
In section \ref{sec:selberg}, which is purely algebraic and independent of the rest of the paper, we compute the Selberg-type integrals on which the proof of Theorem \ref{thm:main_1} relies, and discuss their relation to some known integrals. 
Then in section \ref{sec:general_existence} we use the results of the previous sections, as well as the appendices, to show that every $\OO(p,q)$-invariant valuation is given by an $\OO(p,q)$-invariant Crofton distribution, along the way proving the non-triviality of some Crofton formulas in a few universal families, thus proving Theorems \ref{thm:main_1} and \ref{thm:main_universal_families} .
In section \ref{sec:mu_c} we study the invariant distribution of minimal support, proving the last part of Theorem \ref{thm:main_mu_c}. We rely on results of James and Constantine \cite{constantine}, \cite{constantine_james}, describing the zonal functions on the Grassmannian. Then in section \ref{sec:q=2} we construct explicitly a basis of $\OO(p,2)$-valuations through invariant Crofton distributions, along the way completing the proofs of the $q=2$ statements in Theorems \ref{thm:main_mu_c} and \ref{thm:main_trivial}.
In Appendix \ref{sec:ellipsoid} we compute the projections of the family of ellipsoids used in the proof of Theorem \ref{thm:main_1}.
Finally in Appendix \ref{sec:functorial}, coauthored by T. Wannerer, we study the functorial properties of Crofton distributions, which are used throughout the paper.

\subsection*{Acknowledgments}
I wish to thank Andreas Bernig, for many fruitful discussions on this and related topics; Thomas Wannerer, with whom I was discussing the functorial properties of Crofton measures while enjoying his hospitality in Jena; Semyon Alesker, for introducing me to many of the central notions and ideas in this paper; and Jacob Tsimerman, who made the very helpful suggestion of applying the DiPippo-Howe method for the computation of the integral in section \ref{sec:selberg}. I am also grateful to Andreas, Semyon and Thomas for their numerous helpful comments on the paper. I also would like to thank Peter J. Forrester and Sasha Sodin for useful comments and references on Selberg's integral, Gautier Berck, from whom I learnt about the connection between centro-affine surface area and Crofton formulas, and Dror Bar-Natan for his help with the integral-computing code.
 
\section{Preliminaries}\label{sec:preliminaries}

We will be making use of the Iverson notation $[S]\in\{0,1\}$, which simply equals $1$ if statement $S$ is true and $0$ if it is false. A half-integer is any number $x\in \R$ with $2x\in\mathbb Z$. The half-integers $x\not\in \mathbb Z$ will be called strict half-integers. $V$ will denote a real linear space of dimension $n$, $\Gr_k(V)$ the Grassmannian of non-oriented $k$-dimensional linear subspaces in $V$, and $\AGr_k(V)$ the corresponding Grassmannian of affine subspaces. We reserve the letter $E\in \Gr_k(V)$ to denote a general point in the linear Grassmannian when describing the fiber of a vector bundle. We denote the compact convex subsets of $V$ by $\mathcal K(V)$. The space of Lebesgue measures (densities) on $V$ will be denoted $\Dens(V)$. For a group $G\subset \GL(V)$, $\overline G:=G \ltimes V$ is the group of affine maps generated by $G$ and translations in $V$.
We let $\mathcal M^\infty(X)$,  $\mathcal M(X)$ and $\mathcal M^{-\infty}(X)$ denote the spaces of smooth measures, Borel measures and distributions (generalized measures) on the manifold $X$, respectively. Those are the sections of respective regularity of the bundle of densities $|\omega_X|$ over $X$ that has fiber $\Dens(T_xX)$ over $x\in X$. More generally, given a vector bundle $\mathcal E$ over $X$, the $\mathcal E$-valued measures on $X$ of given regularity are the corresponding spaces of sections: $\mathcal M^\nu(X,\mathcal E)=C^\nu(X, |\omega_X|\otimes \mathcal E)$, $\nu\in\{+\infty, -\infty \}$.

\subsection{Valuation Theory}

In this note we are only concerned with even valuations, which allows to present shorter definitions than in the general case. As the results in this paper are largely independent of the general theory of valuations, we will use  definitions and descriptions that are most easily applicable for our purposes, sometimes masking deep theorems that lie beneath those descriptions. For a survey of the modern theory of valuations, see \cite{alesker_survey}, \cite{bernig_survey}, \cite{fu_survey} and the references therein. 

\begin{Definition}
	 A valuation is a function $\phi:\mathcal K(V)\to \R$ which is a finitely additive measure on convex bodies: $\phi(K\cup L)+\phi(K\cap L)=\phi(K)+\phi(L)$ whenever $K,L,K\cup L\in\mathcal K(V)$.
 	 The space of translation-invariant valuations continuous with respect to the Hausdorff metric on $\mathcal K(V)$ is denoted $\Val(V)$.
 \end{Definition}
Fixing any Euclidean ball $B\subset V$, we get a Banach norm $\|\phi\|=\sup_{K\subset B}||\phi(K)|$. 
One readily decomposes $\Val(V)=\Val^+(V)\oplus \Val^-(V)$, where $\Val^\pm(V)=\{\phi: \phi(-K)=\pm \phi (K)\}$. Denote the $k$-homogeneous valuations by $\Val_k(V)=\{\phi: \phi(\lambda K)=\lambda^k\phi(K), \forall \lambda>0 \}$.
It is easy to see that $\Val_0(V)=\Span\{\chi\}$, where $\chi(K)=1$ is the Euler characteristic.
For $1\leq k\leq n-1$, $\Val_k^{\pm}(V)$ is infinite-dimensional. We recall two classical facts.

\begin{Theorem*}[Hadwiger \cite{hadwiger},\cite{schneider_book14}]
	$ \Val_n(\R^n)=\Span\{\vol_n\}$.
	\end{Theorem*}
		
\begin{Theorem*}[McMullen's decomposition \cite{mcmullen77}]
	$\Val(\R^n)=\oplus_{k=0}^n \Val_k(\R^n)$.
\end{Theorem*}

It is a consequence of Alesker's irreducibility theorem \cite{alesker_mcmullen} that the following definition is equivalent to his original definition of a smooth valuation.

\begin{Definition}
	The space of smooth, translation invariant, even, $k$-homogeneous valuations on $V$, denoted $\Val_k^{+,\infty}(V)$, is the image of the Crofton map 
	\[\Cr: \mathcal M^{\infty}(\AGr_{n-k}(V))^{tr}\to \Val_k^+(V)\]
	given by
	\[\Cr(\mu)(K)=\int_{{\AGr}_{n-k}(V)}\chi(K\cap \overline E)d\mu(\overline E)  \] 
	We say $\mu$ is a smooth Crofton measure for $\phi=\Cr(\mu)$.
	$\Val_k^{+,\infty}(V)$ is equipped with a natural topology that makes it into a Frechet space, and the Crofton map $ \Cr: \mathcal M^{\infty}(\AGr_{n-k}(V))^{tr}\to \Val_k^{+,\infty}(V)$ a continuous surjection.
\end{Definition}
Throughout the paper we use interchangeably the isomorphic spaces 
\[\mathcal M^\infty(\AGr_{n-k}(V))^{tr} = \mathcal M^\infty(\Gr_{n-k}(V), \Dens(V/E))   \]
where the latter is the space of smooth measures with values in the line bundle of densities in the space $V/E$ over $E\in \Gr_{n-k}(V)$. The same applies also to spaces of distributions.
Given a convex body $K\subset V$, this corresponds to the equivalence of Crofton and Kubota formulas:
\[\int_{\AGr_{n-k}(V)}\chi(K\cap \overline E)d\mu(\overline E) = \int_{\Gr_{n-k}(V)} \langle  \pi_{V/E}(K), d\mu(E)\rangle \]

\begin{Definition}
	The Klain map $\Kl: \Val_k^{+,\infty}(V)\to C^\infty (\Gr_k(V), \Dens(E))$ is given by $\Kl(\phi)(E)=\phi|_E$.
	By Klain's theorem \cite{klain}, it is injective.
\end{Definition}

The composition \[T_k=\Kl\circ\Cr :\mathcal M^\infty(\Gr_{n-k}(V), \Dens(V/E)) \to C^\infty (\Gr_k(V), \Dens(E))\] is the well-known cosine transform written in $\GL(V)$-equivariant form. It assumes the more familiar form 
$T_k : C^\infty(\Gr_k(V))\to C^\infty(\Gr_k(V))$ \[T_k(f)(E)=\int_{\Gr_k(V)}f(F)\langle E, F\rangle dF\]
if one fixes a Euclidean structure on $V$ and uses it to identify $E$ with $E^\perp$, and to trivialize both line bundles.

\begin{Definition}
	The space of translation-invariant, even, $k$-homogeneous, generalized valuations $\Val_k^{+,-\infty}(V)$ is defined to be the twisted topological dual $\Val_{n-k}^{+,\infty}(V)^*\otimes \Dens(V) $, equipped with the weak topology.
\end{Definition}

By the Alesker-Poincar\'{e} duality \cite{alesker04_product}, there are natural dense inclusions \[\Val_k^{+,\infty}(V)\subset\Val_k^+(V) \subset\Val_k^{+,-\infty}(V).\]

Generalized valuations can be thought of as valuations on convex bodies with smooth support function, denoted $\mathcal K_s(V)$. Given $K\in\mathcal K_s(V)$, it induces a natural functional on $\Val^{+,-\infty}(V)$, extending the evaluation at $K$ map $\ev_K:\Val(V)\to \R$.

One has an extension of the Klain and Crofton maps:

\[\mathcal M^{-\infty}(\Gr_{n-k}(V), \Dens(V/E)) \mathop{\longrightarrow}^{\mathrm{\Cr}}  \Val^{+,-\infty}_k(V)\mathop{\longrightarrow}^{\mathrm{\Kl}} C^{-\infty}  (\Gr_k(V), \Dens(E))
\]
which are still surjective and injective, respectively. Moreover, this generalized Crofton map is the dual of the smooth Klain map and vice versa, see \cite{alesker_faifman} for details.

We refer to the elements of $\mathcal M^{-\infty}(\Gr_{n-k}(V), \Dens(V/E))$ as $k$-homogeneous Crofton distributions.
When no confusion can arise, elements of $\Val^{+,-\infty}(V)$ will simply be referred to as valuations.
When the space of valuations is twisted by some linear space, we let $\Kl$, $\Cr$ act by the identity on the extra factor.

\subsection{Geometry of $\OO(p,q)$}
Let us start by introducing some notation and definitions which we will use throughout the paper, for details we refer to \cite{bernig_faifman_opq}.

Let $V=\R^n$ be equipped with a non-degenerate quadratic form $Q$ of signature $(p,q)$. We will always assume $q\leq p$. We will write interchangeably $V$ or $\R^{p,q}$, as well as $\OO(Q)$ or $\OO(p,q)$ for the corresponding indefinite orthogonal group. When we do explicit computations, we work with the standard Euclidean form $P(x)=\sum x_j^2$ and $(p,q)$-form $Q(x)=\sum_{j=1}^px_j^2-\sum_{j=p+1}^{n}x_j^2$. Let $\vol_Q\in\Dens(V)$ denote the Lebesgue measure induced by $Q$, which will often be implicitly used to trivialize $\Dens(V)$.

A \textit{$Q$-compatible Euclidean form $P$} is a Euclidean form such that $\OO(P)\cap \OO(Q)$ is a maximal compact subgroup of $\OO(Q)$, and also $\sup_{x\neq 0} \frac{Q(x)}{P(x)}=1$, $\inf_{x\neq 0} \frac{Q(x)}{P(x)}=-1$. One then has an induced $P$- and $Q$- orthogonal decomposition $V=V^p_+\oplus V^q_-$ with $V^\bullet_\pm=\{x:Q(x)=\pm P(x)\}$, $\dim V^p_+=p$ and $\dim V^q_-=q$. We let $E^P$ and $E^Q$ denote the orthogonal complement with respect to the corresponding quadratic form. The linear involution $S_P\in \OO(Q)\cap \OO(P)$ is defined by $Q(u,v)=P(S_Pu,v)$. It holds for a subspace $E\subset V$ that $(E^P)^Q=(E^Q)^P=S_P(E)$.

For $g\in \GL(V)$ we define $\psi_g:\Gr_k(V)\to (0,\infty)$ by 
\[\psi_g(E)=|\Jac g:E\to gE|^{-2}  \] which is computed with respect to the Euclidean structure defined by $P$. 
Equipping $\Gr_k(V)$ with the Euclidean structure induced by $P$ through the identification $T_E\Gr_k(V)=E^*\otimes E^P$ and the Hilbert-Schmidt norm $\|L:E\to E^P\|_{\textrm{HS}}=\tr(LL^{T})$, it holds that $|\Jac_E g: \Gr_k(V)\to\Gr_k(V)|=\psi_g(E)^{\frac n 2}|\det g |^k$.

{ The orbits of $\Gr_k(V)$ under $\OO(Q)$ are given by \[X^k_{a,b}=\{E: Q|_E\text{ has signature } (a,b)\},\] 
where $\max(0, k-q)\leq a\leq p$,	$\max(0, k-p)\leq b\leq q$. We denote the unique closed orbit corresponding to the minimal pair $(a,b)$ by $X^k_c$.

For $E\in X^k_{a,b}$, set $E_{0}=E\cap E^{Q}$, $\dim E_{0}=r=k-(a+b)$. By Proposition 4.2 in \cite{bernig_faifman_opq}, there is a $\Stab(E)$-equivariant
identification 
\begin{equation}\label{eq:normal_bundle_orbit}
N_{E}X_{a,b}^{k}:=T_E\Gr_k(V)/T_EX^k_{a,b}=\Big(\Sym^{2}E_{0}\Big)^{*}
\end{equation}
	 
}

  	\begin{Definition}\label{def:2.1}
  		Set $N=\min(q, k, n-k)$, $N_q=\max(0,k-q)$ and $N_p=\max(0, k-p)$. For $E\in\Gr_k(V)$, choose a $P$-orthonormal basis $B_E=(u_1,\dots, u_N, v_1,\dots, v_{N_q},\\ w_1,\dots, w_{N_p})$ of $E$ such that $v_j\in V^p_+$ and $w_j\in V^q_-$.
  		Let $\lambda_1(E)\geq\dots\geq \lambda_{N}(E)$ be the eigenvalues of $M_P(E):=(Q(u_i,u_j))_{i,j=1}^N$. 
  		Define also $\theta_j(E)\in [0,\frac \pi 2]$ by $\lambda_j(E)=\cos2\theta_j(E)$.
  	\end{Definition}
  	
  	Recall the probability distribution of the principal angles between two random planes in $\R^n$, see e.g. \cite{james_distribution} for details.
  	\begin{Theorem*}[Principal angles distribution]
  		Let $E\in \Gr_e(\R^n)$ be chosen randomly according to the $\SO(n)$-invariant probability measure on the Grassmannian, and let $F\in \Gr_f(\R^n)$ be fixed. Set $m=\min(e,f)$. Let $\mu_j=\cos^2\theta_j$, $j=1,\dots m$ be the ordered squared cosines of the principal angles between $E$ and $F$, that is, $\mu_1\geq\dots\geq \mu_m$ are the eigenvalues of $L^TL$, where $L$ is the orthogonal projection $L:E\to F$ if $e\leq f$ and $L:F\to E$ otherwise, written with respect to some orthonormal bases of $E, F$. 
  		The probability density is proportional to
  		
   		\[ \prod_{i<j} (\mu_i-\mu_j)\prod_{j=1}^m\mu_j^{\frac{|f-e|-1}{2}} \prod_{j=1}^m(1-\mu_j)^{\frac{|n-(e+f)|-1}{2}} d\mu \]
	\end{Theorem*}
  	 \begin{Remark}
  	 	The constant can be deduced from Selberg's integral \cite{selberg}, namely
  	 	\[\int_{[0,1]^m}  \prod_{i<j} |\mu_i-\mu_j|\prod_{j=1}^m\mu_j^{a} \prod_{j=1}^m(1-\mu_j)^{b} d\mu =\frac{\Gamma_m(a+\frac{m+1} 2)\Gamma_m(b+\frac{m+1} 2)\Gamma_m(1+\frac{m}{2})}{\pi^{m^2/2}\Gamma_m(a+b+m+1)} \]
  	 	where the Gamma function of the cone of $m\times m$ positive symmetric matrices is
  	 	\begin{equation}\label{def:gamma_n}\Gamma_m(s):=\int_{X>0}e^{-\tr(X)}(\det X)^{s-\frac{m+1}{2}}dX=\pi^{\frac{m(m-1)}{4}}\prod_{i=0}^{m-1}\Gamma(s-\frac{i}{2})\end{equation}
  	 	
  	 \end{Remark}
  
   	\begin{Proposition}\label{prop:beta_distribution}
    Set $N=\min(q,k, n-k)$. The angles  $\theta_j(E)$, $1\leq j\leq  N$ are the non-trivial principal angles between $E$ and $V^p_+$.
	The eigenvalues $1\geq \lambda_1(E)\geq\dots\geq\lambda_N(E)\geq-1$ of $M_P(E)$ have probability density proportional to  
	\[\prod_{i<j}(\lambda_i-\lambda_j)\prod_{i=1}^N(1-\lambda_i)^{\frac{|q-k|-1} {2}}\prod_{i=1}^N (1+\lambda_i)^{\frac{|p-k|-1}{2}}d\lambda\]
 
  	\end{Proposition}

  	\proof
  	We work in $V=\R^n$ with all the standard structures. Take $E\in\Gr_k(\R^n)$. Let $N_q, N_p$ and $B_E$ be as in Definition \ref{def:2.1}.
  	Write $u_j=u_j^++u_j^-$ with $u_j^+\in\R^{p,0}$, $u_j^-\in \R^{0,q}$. Note that $P(u_i^+, v_j)=0$ 
    
    Let $L$ be the projection operator from $E$ to $\R^{p,0}$, written with respect to the basis $B_E$ in $E$ and the standard basis in $\R^{p,0}$. Note that $L^TL$ is an $(N, N_q, N_p)$ block-diagonal matrix, with diagonal given by the matrices $(P(u_i^+, u_j^+))_{i,j=1}^N$, $I_{N_q}$ and $0_{N_p}$. Now $Q(u_i,u_j)=Q(u_i^+, u_j^+)+Q(u_i^-, u_j^-)=P(u_i^+, u_j^+)-P(u_i^-, u_j^-)=2P(u_i^+,u_j^+)-\delta_i^j$.
	It follows that 
	\[[M_P(E)]_N=(Q(u_i,u_j))_{i,j=1}^N=2[L^TL]_N-I_N\] where $[\bullet]_N$ denotes the $N$-th principal minor. 
	
	Thus the eigenvalues of $[L^TL]_N$ are $\frac{1+\lambda_j(E)}{2}=\cos^2\theta_j(E)$, $j=1,\dots,N$, and their probability distribution is immediate from the distribution of the principal angles. This concludes the proof.
  
  	\endproof

  	\begin{Lemma}\label{lem:MP_orthogonal}
  		For $E\in\Gr_k(V)$ and $1\leq j\leq N$,	$\theta_j(E^Q)=\theta_j(E^P)=\frac \pi 2-\theta_{N+1-j}(E)$.
  	\end{Lemma}
  	\proof
  	The first equality is evident since $E^Q=S_PE^P$ and $S_P\in \OO(P)\cap \OO(Q)$. The second equality follows from Proposition \ref{prop:beta_distribution}.
  	\endproof
 
	\subsection{The classification of $\OO(p,q)$-invariant valuations}\label{sub:opq_classification}
	 We will be relying on the following classification results from \cite{bernig_faifman_opq}.
	 \begin{Theorem}\label{thm:dimension=2}
	 	For all $p\geq q\geq 1$ and $1\leq k\leq p+q-1$, $\dim \Val_k^{-\infty}(\R^{p,q})^{\OO(p,q)}=2$.
	 \end{Theorem}
	 
	 For $a+b=k$, let $\kappa_{a,b}\in C(\Gr_k(\R^{p,q}), \Dens(E))$ be given by $\kappa_{a,b}(E)=\vol_{Q|_E}$ when $E\in X^k_{a,b}$, and zero otherwise. One easily checks that $\kappa_{a,b}$ is in fact continuous.
	 
	 \begin{Theorem}\label{thm:klain}
  	 	For $p\geq q \geq 1$ and $1\leq k\leq p+q-1$, the image of the Klain map on $\Val_k^{-\infty}(\R^{p,q})^{\OO(p,q)}$ is spanned by $\kappa_k^{\textbf{\emph{cos}}}:=\sum_{a+b=k}\cos(\frac\pi 2 b) \kappa_{a, b}$ and $\kappa_k^{\textbf{\emph{sin}}}:=\sum_{a+b=k}\sin(\frac\pi 2 b)\kappa_{a, b}$, where both sums range over all open orbits $X^k_{a,b}$.
 	 \end{Theorem}

	  In \cite{bernig_faifman_opq}, the class of Klain-Schneider valuations was introduced, such that $\Val(V)\subset \Val^{\text{KS}}(V)\subset \Val^{-\infty}(V)$. 
	  For even valuations, those are just the generalized valuations with continuous Klain section. Moreover, the operations of pull-back and push-forward by linear maps were shown to extend to $\Val^{\text{KS}}(V)$. Combined with the classification, one arrives at the following fact.
	  
	 \begin{Theorem}\label{thm:isomorphism}
	  Denote $d=p+q-(p'+q')$. Let $j:\R^{p',q'}\to \R^{p,q}$ be an isometric inclusion, and $\pi:\R^{p,q}\to \R^{p',q'}$ a $Q$-orthogonal projection.  
	  
	  Then 
	  $j^*:\Val_k^{-\infty}(\R^{p,q})^{\OO(p,q)}\to \Val_k^{-\infty}(\R^{p',q'})^{\OO(p',q')}$ and	$\pi_*:\Val_k^{-\infty}(\R^{p,q})^{\OO(p,q)}\\\to \Val_{k-d}^{-\infty}(\R^{p',q'})^{\OO(p',q')}$ are isomorphisms whenever the corresponding target space is $2$-dimensional.
	 \end{Theorem}
	 
	 When $p=q$, we may identify $\R^{p,p}=\C^p$ and take the form $Q(z)=|\Re(z)|^2-|\Im(z)|^2$. Since $i^*Q=-Q$, the pull-back by the imaginary unit 
	 defines an involution of spaces of $\OO(Q)$-invariants. In particular, we obtain an eigenspace decomposition into $\pm1$ eigenspaces of $i^*$:
	 $\Val_k^{-\infty}(\R^{p,p})^{\OO(p,p)}=W_{k,+}^{p,p}\oplus W_{k,-}^{p,p}$, referred to as $i$-even and $i$-odd. Theorem \ref{thm:klain} implies they are both $1$-dimensional.
	 
	 \subsection{Wavefronts of invariant Crofton distributions}
	  For the standard facts on wavefronts, see \cite{duistermaat}, \cite{guillemin_sternberg}.
	  \begin{Proposition}\label{prop:Crofton_opq_wavefront}
	  Let $\mu\in \mathcal M^{-\infty}(\Gr_{k}(V), \Dens(V/E))$ be an $\OO(p,q)$-invariant Crofton distribution. Then
	  \begin{enumerate}
	  	\item $\WF(\mu)\subset\bigcup_{a,b}N^* X^k_{a,b}$.
	  	\item Fix a non-degenerate subspace $U\subset V$ in $V=\R^{p,q}$. Set $S_j=\{E\in\Gr_k(V): \dim E\cap U=j\}$. Then $\WF(\mu)\cap N^*S_j=\emptyset$.
	  \end{enumerate} 
	  \end{Proposition}
	  \proof
	  
	  The first part follows from the $\OO(p,q)$-invariance of $\mu$. For the second part, take $E\in X^k_{a,b}\cap \Gr_{k}(U)$ and set $E_0=E\cap E^Q$. For $T\in T_E\Gr_k(V)$, write $T:E\to V/E$, choose any lift $\tilde T:E\to V$ and decompose $\tilde T=\tilde T_1+\tilde T_2$, $\Im(\tilde T_1)\subset U$ and $\Im(\tilde T_2)\subset U^Q$. Composing with the quotient to $V/E$, we get $T_1\in T_E S_j$ and $T_2\in T_E X^k_{a,b}$, since $U^Q\subset E_0^Q$. Thus $T_E\Gr_k(V)=T_E X^k_{a,b}+T_E S_j$. 
	  	  
	  \endproof
	  
	  \begin{Remark}\label{rem:wavefront_applied}  Combined with Proposition \ref{prop:Crofton_restriction}, we conclude that for an isometric embedding $j:\R^{p',q'}\to\R^{p,q}$ and an $\OO(p,q)$-invariant Crofton distribution $\mu$, the restriction $j^*\mu$ is a well-defined $O(p',q')$-invariant Crofton distribution, and similarly the push-forward $\pi_*\mu$ is well-defined and invariant for a $Q$-orthogonal projection $\pi:\R^{p,q}\to\R^{p',q'}$.
	  \end{Remark}

 \section{The trivial obstruction}\label{sec:trivial_obstruction}
 
 We will make use of the following well-known fact (see Appendix A in \cite{bernig_faifman_opq}).
 Let a smooth Lie group $G$ act on a manifold $X$ with finitely many orbits, all of which are locally closed submanifolds. Let $L$ be a $G$-vector bundle over $X$, and $Y \subset X$ a $G$-invariant locally closed submanifold. For an integer $\alpha\geq 0$, define a new $G$-bundle over $Y$ by  
\begin{displaymath}
F^\alpha_Y|_y=\Sym^\alpha(N_yY) \otimes \Dens^*(N_yY) \otimes L|_y.                                                                                                                                                                                                                                                                                                                                                                                                                                                                                                                                                                                    \end{displaymath}
We will write $\Gamma^{-\infty}_Z(L)$ for generalized sections of $L$ with support in $Z\subset X$.
  \begin{Proposition}
 	\label{prop:LocallyClosedOrbits}

 	Let $Z\subset X$ be a closed $G$-invariant subset.
 	Decompose $Z=\bigcup_{j=1}^J Y_j$ where each $Y_j$ is a $G$-orbit. Then
 	\begin{displaymath}
 	\dim \Gamma^{-\infty}_Z(X,L)^G \leq \sum_{\alpha=0}^\infty \sum _{j=1}^J \dim \Gamma^\infty(Y_j, F^\alpha_{Y_j})^{G}.
 	\end{displaymath}
 	
 	More generally, if $Z_1\subset Z_2$ are two $G$-invariant closed subsets of $X$ then
 	\begin{displaymath}
 	\dim \Gamma^{-\infty}_{Z_2}(X,L)^G \leq \dim \Gamma^{-\infty}_{Z_1}(X,L)^G+\sum_{\alpha=0}^\infty\sum _{Y_j \subset Z_2\setminus Z_1} \dim \Gamma^\infty(Y_j, F^\alpha_{Y_j})^{G}. 
 	\end{displaymath}
 	
\end{Proposition}
 	In particular, fixing $y_j \in Y_j$ we have 
 	\begin{displaymath}
 	\dim \Gamma^{-\infty}(Y_j, F^\alpha_{Y_j}) ^G= \dim \Gamma^{\infty}(Y_j, F^\alpha_{Y_j}) ^G=\dim \left(F^{\alpha}_{Y_j}|_{y_j}\right)^{\Stab(y_j)}.
 	\end{displaymath}

 % % % % % % % % % % % % % %
  A Crofton measure (distribution) for an $(n-k)$-homogeneous even valuation is a (generalized) section of the vector bundle over $\Gr_{k}(V)$ with fiber $\Dens(V/E)\otimes |\omega|_E$, where $|\omega|_E=\Dens(T_E\Gr_{k}(V))$.
  Writing $\Stab(E)\subset \GL(V)$ for the stabilizer, one has the $\Stab(E)$-equivariant identification \[\Dens(V/E)\otimes |\omega|_E=\Dens(V)^{k}\otimes \Dens^*(E)^{n+1}\]
  
  It follows that an $\OO(Q)$-invariant Crofton distribution for $\phi\in \Val_{n-k}^{+,-\infty}(V)$ is given by a  generalized section $\mu\in \Gamma^{-\infty}(\Gr_{k}(V), \Dens^*(E)^{n+1})^{\OO(Q)}$.
  
  Using a Euclidean trivialization of the latter bundle and writing $dE$ for the $\OO(P)$-invariant probability measure on $\Gr_{k}(V)$, we may identify an $\OO(Q)$-invariant Crofton distribution with a generalized function $\hat\mu\in C^{-\infty}(\Gr_{k}(V))$ satisfying $g^*\hat\mu=\psi_g^{-(n+1)/2}\hat{\mu}$. Alternatively, since $g_*(dE)=\psi_{g^{-1}}(E)^{n/2}dE$, we can identify an invariant Crofton distribution with a distribution $\tilde\mu\in \mathcal M^{-\infty}(\Gr_k(V))$ s.t. $g_*\tilde\mu=\psi_{g^{-1}}^{-1/2}\tilde{\mu}$.
 
  Denote $X^k_{\leq m }=\cup_{a+b\leq m}X^k_{a,b}$, which is a compact subset of $\Gr_k(V)$. Similarly define the open submanifold $X^k_{\geq m }$, and the locally closed submanifold $X^k_{=m}$.	Define the space $\Cr(r)$ of all Crofton distributions supported on $X^k_{\leq k-r}$. We can now formulate the main result of the section, which gives an upper bound on the dimension of $\OO(Q)$-invariant Crofton distributions.
 	
 \begin{Proposition}\label{prop:generalCroftonBound}
 	There is a family $\mu_{a,b}\in \mathcal M_{X^k_{a,b}}^{-\infty}(\Gr_k(V)\setminus X^k_{\leq a+b-1}, \Dens(V/E))$ attached to the various orbits, such that $\supp \mu_{a,b}\subset X^k_{a,b}$, with the following properties.
 	\begin{enumerate}
 		\item If $a+b<k$ and $k-(a+b)\equiv n+1\mod 2$ then $\mu_{a,b}=0$.
 		\item For all $s\leq k$, \[\mathcal M^{-\infty}_{ X^k_{\leq s}}\left(\Gr_k(V)\setminus X^k_{\leq s-1}, \Dens(V/E)\right)^{\OO(Q)}=\Span (\mu_{a,b})_{a+b=s}\] 
 	\end{enumerate}
 	
   \end{Proposition}
   \begin{Remark}
	   	For the open orbits, $\mu_{a,b}$ is easy to describe through the Euclidean trivialization: as a generalized function, it equals $|\det M_P(E)|^{-\frac{n+1}{2}}$ for $E\in X^k_{a,b}$ and zero elsewhere. It is not a-priori clear which of the other $\mu_{a,b}$ are non-vanishing.
   \end{Remark}
 \proof
 Take any $a,b$ with $a+b=s<k$, and denote $r=k-s$. Fix $E\in X_{a,b}^{k}$.
 For any integer $\alpha\geq0$ consider the $\Stab(E)$-module
 \[F^{\alpha}_{E}= \Sym^{\alpha}(N_{E}X_{a,b}^{k})\otimes \Dens^{*}(N_{E}X_{a,b}^{k})\otimes\Dens^*(E)^{n+1}.\]
 Let us study the invariants of $F^{\alpha}_{E}$. All other factors being one-dimensional, the existence of an invariant implies the existence
 of an invariant 1-dimensional subspace of $\Sym^{\alpha}(N_{E}X_{a,b}^{k}).$
 Set $E_{0}=E\cap E^{Q}$, $\dim E_{0}=r$. By equation (\ref{eq:normal_bundle_orbit})  we have 
$N_{E}X_{a,b}^{k}=\Big(\Sym^{2}E_{0}\Big)^{*}$. By Witt's extension theorem \cite{witt}, the representation $\Stab(E)\to \GL(E_{0})$
 is onto $\GL(E_{0})$. 
 
 Observe that $\Sym^{\alpha}(N_{E}X_{a,b}^{k})=\Sym^{\alpha}(\Sym^{2}E_{0})^{*}$
 has precisely one 1-dimensional invariant subspace when $\alpha=mr$,
 $m\geq0$, spanned by the $m$-th power of the determinant on $\Sym^2 E_0$, s.t. $g\in \GL(E_{0})$ acts on this subspace by
 $\det(g)^{-2m}$, while for other values of $\alpha$ there are no
 invariant one-dimensional subspaces.  To see this, recall that a one-dimensional representation of $\GL(\R^r)$ 
 factorizes through the determinant, so $\Sym^\alpha \Sym^2 (\R^r)$ can only contain the unique isomorphism class of a one-dimensional representation given by $\det^\frac{2\alpha}{r}$, and at the same time every irreducible summand in $\Sym^\alpha \Sym^2(\R^r)$ appears with multiplicity one (by a result of Thrall \cite{thrall}, see also Howe \cite{howe} page 562).
 
 Note that $F^{mr}|_{E}$ is $\Stab(E)$-isomorphic to 
 \[
 \Dens^{*}(E/E_0)^{n+1}\otimes \Dens^{*}(E_{0})^{n+1}\otimes \Dens(\Sym^{2}E_{0})\otimes \Sym^{mr}(\Sym^{2}E_{0})^{*}
 \]
 There is a non-degenerate $\Stab(E)$-invariant quadratic form
 induced on $E/E_{0}$, and so $\Stab(E)$ acts on
 the first factor trivially. Consider an element of $\Stab(E)$
 acting on $E_{0}$ by the scalar $\lambda>0$. Then its action
 on $\Dens^{*}(E_{0})^{n+1}\otimes\Dens(\Sym^{2}E_{0})\otimes \Sym^{mr}(\Sym^{2}E_{0})^{*}$
 is by the scalar $\lambda^{r(n+1)}\lambda^{-2\frac{r(r+1)}{2}}\lambda^{-2rm}.$
 We conclude that a non-trivial invariant of $F^\alpha_E$ exists if and only if $\alpha=mr$ and $r(n+1)-r(r+1)-2rm=0\iff m=\frac{n-r}2$. 
 
 When $n-r$ is odd, the second part of Proposition \ref{prop:LocallyClosedOrbits} with $Z_1=X^k_{\leq k-r-1}$ and $Z_2= X^k_{=k-r}$ implies that $\Cr(r)^{\OO(Q)}=\Cr(r+1)^{\OO(Q)}$.
 When $n-r$ is even, take $Z_2=\cup_{a+b\leq k-r}X^k_{a,b}$, $Z_1=Z_2\setminus X^k_{a_0,b_0}$, and apply again 
 Proposition \ref{prop:LocallyClosedOrbits}.
 \endproof

 \begin{Remark}\label{rem:order_of_mu}
 	It follows from the proof that any $\mu\in \Cr(r)^{\OO(Q)}$ restricted to $X^k_{\geq k-r}$ is  a distribution of order $\alpha=\frac{r(n-r)}{2}$.
   	\\We also see that when $n\not\equiv \min(k, n-k, q)\mod 2$, there is no $\OO(Q)$-invariant Crofton distribution supported on $X^k_c$.

 \end{Remark}

 \section{Constructing $\OO(p,q)$-invariant Crofton distributions}\label{sec:muro}
  
 Take $N=\min(k,n-k,q)$. Let $P$ be a $Q$-compatible Euclidean structure, and $U\subset \Gr_k(V)$ an open set. Let $B_E=(u_1(E),\dots,u_N(E), v_1(E),\dots, v_{N_q}(E), w_1(E),\\ \dots, w_{N_p}(E)):U\to V^k$ be a smooth field of $P$-orthonormal bases of $E\in U$ as in Definition \ref{def:2.1}. Recall the function $M_P:U\to \Sym_N(\mathbb R)$ given by $M_P(E)=Q(u_i(E),u_j(E))_{i,j=1}^N$, where $\Sym_N(\mathbb R)$ is the space of symmetric $N\times N$ real matrices.  Denote by $U_P\subset\Gr_k(V)$ the open and dense subset of subspaces intersecting both $V^p_+$ and $V^q_-$ generically.

 \begin{Lemma}\label{lem:matrix_is_submersion}
 	$M_P$ is a proper submersion at every $E\in U\cap U_P$.
 \end{Lemma}
 \begin{Remark}
 	It is easy to see that $E\in U_P$ if and only if $M_P(E)$ has no eigenvalue equal to $\pm1$, if and only if $E$ intersects $(E^P)^Q=S_PE$ generically.
  	
 \end{Remark}
 \proof
 Consider a curve $\gamma_1$ through $E$ given by \[\gamma_1(t)=\Span(u_1(t),u_2,\dots,u_N, v_1,..., v_{N_q}, w_1, \dots, w_{N_p})\] with all vectors fixed except for $u_1$, and $\xi=\dot u_1(0)\in E^P$ arbitrary. It follows that 
 \[D_E M_P (\dot{\gamma_1}(0))=
 \begin{pmatrix}
 Q(\xi, 2u_1) & Q(\xi, u_2) & \cdots & Q(\xi,u_N) \\
 Q(\xi,u_2) & 0 & \cdots & 0 \\
 \vdots  & \vdots  & \ddots & \vdots  \\
 Q(\xi,u_N) & 0 & \cdots & 0
 \end{pmatrix}
 \]
 Note that $E\cap(E^P)^Q=E\cap S_PE=\Span(v_1,..., v_{k-q}, w_1, ..., w_{k-p})$ by the generic intersection assumption. Hence $Q(\xi, 2u_1), Q(\xi,u_2),\cdots, Q(\xi,u_N)$ are linearly independent functionals in $\xi\in E^P$, and so the first row of $D_\Lambda M_P(\dot{\gamma_1})$ is arbitrary while the other entries in upper triangle vanish. Replacing  $\gamma_1$ with $\gamma_j$ in the obvious way, we conclude $D_E M_P(\sum \alpha_j \dot\gamma_j(0))$ can be arbitrary, thus concluding the proof.
 \endproof
 
 \begin{Lemma}\label{lem:compatible_cover}
 	One can choose finitely many $Q$-compatible Euclidean structures $P_i$ s.t. $\{U_{P_i}\}$ cover $\Gr_k(V)$.
 \end{Lemma}
 \proof
 Given $E\in \Gr_k(V)$, $M_P$ is submersive at $E$ for a generic choice of $Q$-compatible Euclidean form $P$ by a trivial dimension count. The claim then follows by the compactness of $\Gr_k(V)$.
 \endproof

 The space $\Sym_N(\mathbb R)$ is acted upon by $\GL(N)$ through $g(X)=gXg^T$.
 The following fact is well-known, see e.g. \cite{muro} or \cite{blind}, and goes back essentially to Cayley \cite{cayley} and G$\mathring{\text a}$rding \cite{garding}, who adapted Cayley's formula to the symmetric case. From a modern perspective, it is an instance of the Bernstein-Sato theorem on the meromorphic extension of $|P|^s$. Let $S_{a,b}\subset\Sym_N(\R)$ denote the matrices of signature $(a,b)$. Note that $\overline{S_{a,b}}=\bigcup_{a'\leq a,b'\leq b}S_{a',b'}$.
 
 \begin{Theorem}[G$\mathring{\text a}$rding]
 	For any two integers $a,b\geq 0$ with $a+b=N$, there is a meromorphic in $s$ family of generalized functions $\Phi_a(s)\in C^{-\infty}(\Sym_N(\mathbb R))$ supported on $\overline {S_{a,b}}$, and restricting to $\Phi_a(s)(X)=|\det X|^s$ on $S_{a,b}$. The poles occur precisely at the half-integers $s\leq -1$.
 \end{Theorem}
 
 The following characterization of $\Phi_a(s)$ is a particular case of Muro's Theorem 5.6 in \cite{muro_prehomogeneous}.
 
 \begin{Theorem}[Muro]\label{thm:muro_classification}
 	For any $s_0\in \mathbb C$, the space of generalized functions $D_{s_0}=\{f\in C^{-\infty}(\Sym_N(\R)): g^*f=|\det g|^{2s_0}f, \forall g\in \GL(N)\}$ is of dimension $N+1$. For $s_0\notin \frac{\mathbb Z}{2}\cap(-\infty -1]$, $D_{s_0}$ is spanned by $\Phi_a(s_0)$. At a half-integer $s_0\leq -1$, it is spanned by the leading Laurent coefficients of the various linear combinations of $(\Phi_{a}(s))_{a=0}^N$ at $s=s_0$.
 \end{Theorem}
 
 \begin{Definition}
 	Denote by $S^r\subset Sym_N(\mathbb R)$ the collection of symmetric matrices of rank $N-r$, and by $S^r_\pm\subset S^r$ the subsets of positive/negative semi-definite matrices.
 	Denote by $w_a(s_0)$ the order of the pole of $\Phi_a(s)$ at $s_0$, and by $\Psi_a(s_0)$ the corresponding leading Laurent coefficient.
 \end{Definition} 
   A description of the poles and Laurent coefficients of $\Phi_a$ and their linear combinations was obtained by Muro in \cite{muro} using Sato's hyperfunctions, and later by Blind \cite{blind} using microlocal methods. Here is what we will need.
 
 \begin{Theorem}[Muro]\label{thm:muro}
 	\begin{enumerate}
		\item $|\det X|^s=\sum_{a=0}^N\Phi_a(s)(X)$ is analytic at even integers and has a simple pole at odd integers $s\leq-1$, while $\sign(\det X)|\det X|^s=\sum_{a=0}^N (-1)^{N-a} \Phi_a(s)$ is analytic at odd integers and has a simple pole at even integers $s\leq -2$. The supports of the residues have positive codimension.
		\item The linear combinations $\sum_{j=0}^{N/2} (-1)^j\Phi_{N-2j}(s)$, $\sum_{j=0}^{N/2}(-1)^j\Phi_{N-1-2j}(s)$ are analytic at $s\notin\mathbb Z$.
		\item  For $s_0=-\frac{n+1}{2}$ with $n>N$, $w_N(s_0)=\lfloor \frac {N+[s_0\in\mathbb Z]}{2}\rfloor$.
		\item  It holds that $\supp(\Psi_a(s_0))\subset \overline{S^{2w_a(s_0)-[s_0\in\mathbb Z]}}$. 
		\item For $N=2$ and $s_0\in\mathbb Z$ we have $\supp\Psi_2(s_0)=S^1_+\cup S^0$,  $\supp\Psi_0(s_0)=S^1_-\cup S^0$,  $\supp\Psi_1(s_0)=S^1\cup S^0$. 
	 
	\end{enumerate}
		
  \end{Theorem}
We now explain how to pull-back $\Phi_a(s)$ using the locally-defined submersion $M_P$ to obtain some $\OO(Q)$-invariant Crofton distributions.
 \begin{Definition}
 	For $s\in\mathbb C$, let $\mathcal D^s$ be the line bundle of $s$-densities over $\Gr_k(V)$, which has fiber $\Dens^s(E)$ over $E\in \Gr_k(V)$. 
 	We say that a choice of section $f(s)\in\Gamma^{-\infty}(U, \mathcal D^s)$ over $U\subset \Gr_k(V)$ for $s\in\Omega\subset\C$ is meromorphic in $s$ if, having fixed a Euclidean metric $P$ and using it to identify all bundles $\mathcal D^s$, one obtains a map $f_P:\Omega\to  C^{-\infty}(U)$ which is meromorphic in $s$. 
 	
 	We denote by $\mathfrak{M}^{-\infty}(\mathcal D^s)$ the sheaf for which $\Gamma(U, \mathfrak{M}^{-\infty}(\mathcal D^s))$ is the space of meromorphic in $s$ maps $\mathbb C\to \Gamma^{-\infty}(U, \mathcal D^s)$. 
 	
 \end{Definition}

 \begin{Proposition}\label{prop:muro_pullback}
 	For every pair of non-negative integers $(a,b)$ with $a+b=k$, $a\leq p$, $b\leq q$, there is a global section $f_a=M_P^*\Phi_{a-\max(0,k-q)}(s)$ of $\mathfrak{M}^{-\infty}(\mathcal D^s)$ supported on $\overline {X^k_{a,b}}$, s.t. whenever $s$ is not a pole of $f_a$,  $f_a(s)$ is $\OO(p,q)$-invariant.
 \end{Proposition}
 \proof
 Let $P_i$ be a collection of $Q$-compatible Euclidean structures as in Lemma \ref{lem:compatible_cover}, inducing the decompositions $V=V_{p,i}\oplus V_{q,i}$, and denote by $U_i=U_{P_i}\subset\Gr_k(V)$ the subspaces intersecting $V_{p,i}\cup V_{q, i}$ generically. For each $i$, cover $U_i$ by open sets $U_{ij}\subset U_i$ so that $M_{ij}=M_{P_i}:U_{ij}\to \Sym_N(\mathbb R)$ can be defined by some smooth field of orthonormal bases of $E$ over $E\in U_{ij}$. Now since $M_{ij}$ is a proper submersion, one obtains a meromorphic in $s$ family of functions $\tilde f_{ij}(E; a, s)\in C^{-\infty} (U_{ij})$ given by $\tilde f_{ij}(\bullet; a, s)=M_{ij}^*\Phi_{a_0}(s)$, where $a_0=a$ if $k\leq q$, and $a_0=a-(k-q)$ if $k\geq q$.
 
 It then obviously holds that on $U_{ij}\cap U_{ij'}$, $\tilde f_{ij}(\Lambda; a,s)$ and $\tilde f_{ij'}(\Lambda; a,s)$ coincide as continuous functions for $\Re(s)>0$. Therefore, they coincide on $U_{ij}\cap U_{ij'}$ as meromorphic functions, and we may merge all $\tilde f_{ij}$ into one meromorphic family $\tilde f_i(\bullet; a,s)\in C^{-\infty} (U_i)$. The corresponding (through $P_i$) section $f_i\in \Gamma(U_i, \mathfrak{M}^{-\infty}(D^s))$ is obviously $\OO(Q)\cap \OO(P_i)$-invariant. Moreover, it is $\mathfrak{so}(Q)$-invariant.

 Next, we claim that $f_i$ and $f_{i'}$ coincide on $U_i\cap U_{i'}$. Since both are meromorphic, we may assume in the following that $\Re(s)>0$
 
 It is easy to see, using Proposition \ref{prop:LocallyClosedOrbits} as in the proof of Proposition \ref{prop:generalCroftonBound}, that for $\Re(s)>0$, no $\OO(Q)$-invariant generalized sections of $\mathcal D^s$ can be supported on a set of positive codimension. It follows that the space of $\OO(Q)$-invariants in $\Gamma^{-\infty}(\Gr_k(V), \mathcal D^s)$ supported on ${\overline{X^k_{a,b}}}$ is at most 1-dimensional.

 Since $U_i\subset \Gr_k(V)$ is dense, it follows by construction that for $\Re(s)>0$, ${f_i}(\bullet; a, s)$ extends by continuity to an $\OO(Q)$-invariant section of $\mathcal D^s$ over $\Gr_k(V)$ supported on $\overline{X^k_{a,b}}$, and by the previous paragraph we can find meromorphic functions $c_i(s)$, such that $c_1(s)=1$ and the sections $c_i(s) f_i(\Lambda; a, s)$ all coincide. Denoting by $p_i(s)\in\Gamma^\infty (\Gr_k(V), \mathcal D^s)$ the Euclidean section defined by $P_i$,  it holds that \[f_i(E;a,s)= |M_{P_i}(E)|^{s/2}p_i(s)\] for $E\in X^k_{a,b}$, so that \[|M_{P_1}(E)|^{s/2}p_1(s)=c_i(s) |\det M_{P_i}(E)|^{s/2}p_i(s)\]
 implying 
 \[c_i(s)=\left(\frac{|\det M_{P_1}(E)|^{1/2}}{|\det M_{P_i}(E)|^{1/2}}\frac{p_1(1)}{p_i(1)}\right)^s\] 
 for all $E\in X^k_{a,b}$. Since $c_i(s)$ is independent of $E$, one has $c_i(s)=c_i^s$ for some $c_i>0$. 
 
 Finally, for $s=1$, $\mathcal D^1$ is the bundle of Lebesgue measures, and it is easy to see that all the extensions of $f_i$ coincide: Assume for simplicity $(a,b)=(k,0)$, take $B_1\subset E \subset V_{p,1}$ a $P_1$-unit cube, and choose $g\in \OO(Q)$ with $g(V_{p,1})=V_{p,i}$. Then $g(B_1)\subset gE$ is a $P_i$-unit cube, and therefore $f_i(1) (gB_1)=1=f_1(1)(B_1)$, but $f_1(1)(B_1)=f_1(1)(gB_1)$ by $\OO(Q)$-invariance of $f_1(1)$. It follows $c_i=1\Rightarrow c_i(s)\equiv 1$. Thus we have shown that $f_i$, $f_{i'}$ coincide in $\Gamma(U_i\cap U_{i'}, \mathfrak M^{-\infty}(\mathcal D^s))$.
 
 We conclude that one has a globally defined section $f_a$ of $\mathfrak{M}^{-\infty}(\mathcal D^s)$ which is $\mathfrak{so}(Q)$-invariant and supported on $\overline {X^k_{a,b}}$. For $\OO(Q)$-invariance, we observe it holds for $\Re(s)>0$ and then invoke uniqueness of meromorphic continuation.
 
 \endproof
 Set $N=\min(q,k,n-k)$, $s_0=-\frac{n+1}{2}$, $a_m=\max(0, k-q)$, $a_M=\min(k, p)$.  
  
 \begin{Definition}\label{def:epsilons}
 	For $\nu\in\{\textbf{\emph {abs}}, \textbf{\emph {sgn}}, \textbf{\emph {cos}}, \textbf{\emph {sin}} \}$, define
 	the functions $\epsilon_\nu:\mathbb Z\to \mathbb R$ by $\epsilon_\textbf{\emph {abs}}(b)=1$, $\epsilon_\textbf{\emph {sgn}}(b)=(-1)^b$, $\epsilon_\textbf{\emph {cos}}(b)=\cos(\frac \pi 2 b)$,  $\epsilon_\textbf{\emph {sin}}(b)=\sin(\frac \pi 2 b)$.
 	\\
 	Define the linear combinations
 	\begin{equation}\label{eqn:mu_01}\mu_{\nu}(s)=\sum_{a_m\leq k-b\leq a_M}\epsilon_\nu(b)f_{k-b}(s)\end{equation}
 	When we need to emphasize the degree of homogeneity of the Crofton distribution (which is the degree of the corresponding valuation), we will write $\mu_\nu^{n-k}(s)$.
 \end{Definition}
 
 It is easy to describe the image of those Crofton distributions under the $Q$-orthogonal complement, denoted here by $\mathbb F_Q$.
 \begin{Lemma}\label{lem:Crofton_orthogonal_relations} The orthogonal complement acts as follows.
 	\begin{itemize}
	 	\item $\mathbb F_Q \mu^{n-k}_\textbf{\emph{abs}}(s)=\mu^{k}_\textbf{\emph{abs}}(s)$.
	 	\item $\mathbb F_Q \mu^{n-k}_\textbf{\emph{sgn}}(s)=(-1)^q\mu^{k}_\textbf{\emph{sgn}}(s)$.
	 	\item $\mathbb F_Q (\mu^{n-k}_\textbf{\emph {cos}}(s)+\sqrt{-1}\mu^{n-k}_\textbf{\emph {sin}}(s))=\sqrt{-1}^{q} (\mu^k_\textbf{\emph {cos}}(s)-\sqrt{-1}\mu^k_\textbf{\emph {sin}}(s))$.
 	\end{itemize}
 \end{Lemma}
 \proof
 Immediate from the definitions together with Lemma \ref{lem:MP_orthogonal}, which implies that $\mathbb F_Q (M_P^*\Phi_a(s))=M_P^*\Phi_{N-a}(s)$.
 \endproof
We can now prove our first main result.  
 \begin{Theorem}\label{thm:1}
 Fix $1\leq k\leq n-1$, $N=\min(q,k, n-k)$, $s_0=-\frac{n+1}{2}$. The dimension of the space of $\OO(Q)$-invariant Crofton distributons is bounded from below by the number of open orbits, namely $N+1$, and the corresponding distributions are spanned by the leading Laurent coefficients at $s_0$ of the linear combinations of $f_a(s)$, $a_m\leq a\leq a_M$. Moreover,
   \begin{enumerate}
   	\item If $n\equiv 3\mod 4$, $\mu_\textbf{\emph{abs}}(s)$ is analytic at $s_0$, while $\mu_\textbf{\emph{sgn}}(s)$ has a simple pole at $s_0$.  Both $\mu_\textbf{\emph{abs}}(s_0)$ and $\Res_{s_0}\mu_\textbf{\emph{sgn}}(s)$ define $\OO(Q)$-invariant Crofton distributions, and they are linearly independent.
	\item If $n\equiv 1\mod 4$, $\mu_\textbf{\emph{sgn}}(s)$ is analytic at $s_0$, while $\mu_\textbf{\emph{abs}}(s)$ has a simple pole at $s_0$.  Both $\mu_\textbf{\emph{sgn}}(s_0)$ and $\Res_{s_0}\mu_\textbf{\emph{abs}}(s)$ define $\OO(Q)$-invariant Crofton distributions, and they are linearly independent.
  	\item If $n\equiv 0\mod 2$, $\mu_{\textbf {\emph{cos}}}(s)$ and	$\mu_{\textbf {\emph{sin}}}(s)$ are analytic at $s_0$ and define there $\OO(Q)$-invariant Crofton distributions, and they are linearly independent.
   	\item If $n\equiv N\mod 2$, there is a one-dimensional space of $\OO(Q)$-invariant Crofton distributions supported on $X^k_c$, spanned by $\mu_c$ which is the leading Laurent coefficient of $f_{a_M}(s)$ at $s_0$. If $n\not\equiv N\mod 2$, there is no such non-trivial invariant Crofton distribution.

  \end{enumerate}

 \end{Theorem}

 \begin{Remark}
 	 It is likely that one can apply the methods of Ricci and Stein \cite{ricci_stein} and Muro \cite{muro_prehomogeneous} to prove that in fact the number of invariant Crofton distributions equals the number of open orbits. This is certainly true for $q=1$, and we also prove that when $q=2$ and $p$ is even, see Proposition \ref{prop:crofton_classification}.
 \end{Remark}

  \proof
 Note that the pull-back by $M_P$ is injective, so that the orders of poles of $f_a(s)$ match those of $\Phi_{a-a_m}(s)$. 
 The first part now follows from Theorem \ref{thm:muro_classification} and Proposition \ref{prop:muro_pullback}.
 Statements i)-iii) then follow immediately from Theorem \ref{thm:muro} and Proposition \ref{prop:muro_pullback}. The linear independence in the first three statements follows from examining the supports: the analytic extensions have support with non-empty interior, while every residue is supported on a subset of positive codimension. In iii), the supports of $\mu_{\textbf {{cos}}}$ and $\mu_{\textbf {{sin}}}$ have each non-empty interior but intersect at a subset of positive codimension.  
 Finally for statement iv), we use Proposition \ref{prop:generalCroftonBound} and items iii) and iv) of Theorem \ref{thm:muro} to define $\mu_c:=M_P^*\Psi_{a_M-a_m}(s_0)$.

 \endproof
 
 \section{A Selberg-type integral}\label{sec:selberg}
 
 Recall that $S_{a,b}\subset\Sym_n(\R)$ denotes the symmetric matrices of signature $(a,b)$. Denote $S_{a,b}(1):=\{X\in S_{a,b}: -I_n\leq X\leq I_n\}$. Due to the multitude of indices in this section, we will write $\sqrt{-1}$ for the imaginary unit.
 	For a function $\epsilon:\mathbb Z\to \C$, define \[D^\epsilon_n(s):=\sum_{a+b=n}\epsilon(b)\int_{S_{a,b}(1)}|\det X|^s dX.\]  $D^\epsilon_n(s)$ is known to be a rational function of $s$, given by the integral above when $\Re(s)$ is sufficiently large.  We will compute $D_n^\epsilon(s)$ explicitly for the four functions  $\epsilon_\textbf{abs}$, $\epsilon_\textbf{sgn}$, $\epsilon_\textbf{cos}$,  $\epsilon_\textbf{sin}$ from Definition \ref{def:epsilons}.
 	
 	The integral over the positive-definite matrices \[D_n^+(s)=\int_{0\leq X\leq I_n}|\det X|^s dX\] and its many generalizations have been considered, sometimes independently, by different authors. We mention two - Selberg \cite{selberg} and  Constantine \cite{constantine}; see also the survey by Forrester and Warnaar \cite{forrester} for an exhaustive overview of this subject.  The value of $D_n^+(s)$ is given by 
 	\[D_n^+(s)=\frac{\Gamma_n(s+\frac{n+1}{2})\Gamma_n(\frac{n+1}{2})}{\Gamma_n(s+n+1)} \]
 	where $\Gamma_n(x)=\pi^{\frac{n(n-1)}{4}}\prod_{i=0}^{n-1}\Gamma(x-\frac{i}{2})$.
 	
 	It appears however that this and similar integrals with arbitrary signatures replacing the positive-definite one have not been considered earlier. Despite the superficial similarity, even the computation of $D_n^{\textbf{abs}}(s)$ is not a straightforward reduction to the positive-definite case as one might expect. For example, $D_n^+(s)$ can be considered as a particular case of Selberg's integral and computed as such, but it is not clear whether any $D_n^{\epsilon}(s)$ fits into a similar, explicitly computable family of integrals.
 	
 	Nevertheless, $D_n^+(s)$ fits also into a different family of integrals, which was considered by Robbins \cite{robbins} and DiPippo-Howe \cite{dipippo_howe}. Denote $\Delta_n^+=\{1\geq x_1\geq \dots\geq x_n\geq 0\}$.
 	Let us quote their result.  
 	\begin{Proposition}[Robbins, DiPippo-Howe]\label{prop:dipippo_howe}
 		For positive real numbers $e_1,\dots,e_n$ one has
 		\[\int_{\Delta_n^+} \det(x_i^{e_j-1})dx=\frac{1}{e_1\cdots e_n}\prod_{1\leq i<j\leq n}\frac{e_i-e_j}{e_i+e_j} \]
 	\end{Proposition}
 	Note that taking $e_j=s+j$, we get a Vandermonde determinant: $\det(x_i^{s+j-1})= \prod_{i<j}(x_j-x_i) \prod_{i=1}^n x_i^s$. This readily allows to compute $D_n^+(s)$.
 	
 	Here we analyze a related family of integrals, that includes $D^\epsilon_n(s)$ at different infinite sequences of values of $s$ depending on $\epsilon$. This is then sufficient to determine the value of $D^\epsilon_n(s)$ at all $s$.  Denote $\Delta_n=\{1\geq x_1\geq\dots\geq x_n\geq -1\}$ and  define the family of integrals \[f_n(e_1,\dots,e_n)=\int_{\Delta_n} \det (x_i^{e_j-1})dx\]
 	where $e_j\geq \frac12$ are either all integers or all strict half-integers (in the latter case, the precise definition appears below). One immediately observes that for a permutation $\sigma\in S_n$, $f_n(\sigma e)=(-1)^\sigma f_n(e)$. We now restate Theorem \ref{thm:integral} in a form better adapted for the proof. 
 	\begin{Proposition}\label{prop:A5}
 		Take $e=(e_1,\dots, e_n)\in \mathbb N^n$. Let $N_+(e)$ be the number of even entries of $e$, and $N_-(e)$ the number of odd entries. 
 		\begin{enumerate}
 			\item If $N_-(e)-N_+(e)\notin\{0,1\}$ then $f_n(e)=0$.
 			\item Assume $N_-(e)-N_+(e)\in\{0,1\}$. Denote $m=\lceil \frac n 2\rceil$, and assume the entries of $e$ are arranged such that $e_i$ is odd for $i\leq m$ and even for $i>m$. Then

 			\begin{equation}\label{eq:super_awesome_closed_form_integral} f_n(e)= \epsilon_n2^n\frac {\displaystyle \prod_{1\leq i<j\leq m}(e_i-e_j) \prod_{m<k< l\leq n} (e_k-e_l)}{\displaystyle\prod_{1\leq i\leq m} e_i \prod_{\scriptsize\begin{array}{c}1\leq i\leq m\\m<k\leq n\end{array}} (e_i+e_k)}
 			\end{equation}
 			where $\epsilon_n\in\{-1, 1\}$ is an 8-periodic sequence given by $\epsilon_1=1$, $\epsilon_{2m+1}=\epsilon_{2m}=(-1)^m\epsilon_{2m-1}$.
 		\end{enumerate}
 	\end{Proposition} 
 	
 	Next we consider $f_n(e_1,\dots, e_n)$ with $e_j>0$ strict half-integers. We use the convention that for $x<0$, $x^{\frac{2k+1}{2}}:=(-1)^k|x|^{\frac{2k+1}{2}}\sqrt{-1}$. 
 	
 	\begin{Proposition}\label{prop:mixed_signs}
 		Let $e=(e_1,\dots, e_n)$ be a vector with positive, strict half-integer coordinates, and let $0\leq m\leq n$ be arbitrary. Assume that $2e_j\equiv 1\mod 4$ for $1\leq j \leq m$ and $2e_j\equiv 3\mod 4$ for $m+1\leq j\leq n$.
 		Then
 		\[f_n(e)=\frac{\delta_n(m)}{\prod_{i=1}^n e_i}\prod_{\scriptsize\begin{array}{c} i<j \\ 2e_i\equiv 2e_j\mod 4\end{array}}\frac{e_i-e_j}{e_i+e_j} \]
 		where \[\delta_n(m)=2^{\frac n 2}\left\{\begin{array}{cc} 1,& n\equiv m\equiv 0\mod 2\\\sqrt {-1},& n\equiv 0, m\equiv 1\mod 2\\e^{(-1)^m\frac {\pi}{4}\sqrt {-1}},& n\equiv 1\mod 2
 		\end{array}  \right. \]
 	\end{Proposition}

 	Let us first see how these result can be used to compute the various $D^{\epsilon}_n(s)$. First we consider the corollaries of Proposition \ref{prop:A5}.
 	\begin{Proposition}\label{prop:Beta_integral_matrix} 
 		\[D_n^{\textbf{\emph {abs}}}(s)=\frac{n!\pi^{\frac{n^2}{2}}}{\Gamma_n(\frac{n+2}2)} \frac{ \displaystyle\prod_{\scriptsize\begin{array}{c}1\leq i< j\leq n\\i\equiv j\mod 2 \end{array} }(j-i) }   {\displaystyle\prod_{\scriptsize\begin{array}{c} 1\leq i\leq n\\i\equiv 1\mod 2 \end{array}} (s+i) \prod_{\scriptsize\begin{array}{c}1\leq i< j\leq n\\i\not\equiv j\mod 2 \end{array}} (2s+i+j)}\]
 		
 		\[D^{\textbf{\emph {sgn}}}_{2m}(s) = \frac{(2m)!\pi^{2m^2}}{\Gamma_{2m}(m+1)}(-1)^m 2^{2m} \frac{ \displaystyle\prod_{\scriptsize\begin{array}{c}1\leq i< j\leq 2m\\i\equiv j\mod 2 \end{array} }(j-i) }   {\displaystyle\prod_{i=1}^m (s+2i) \prod_{\scriptsize\begin{array}{c}1\leq i< j\leq 2m\\i\not\equiv j\mod 2 \end{array}} (2s+i+j)}  \]
 	\end{Proposition}
 	
 	\proof
 	Let $\lambda:\Sym_n(\R)\to\{\lambda_1\geq\dots\geq \lambda_n\}\subset \R^n$ be the spectrum map, mapping a matrix to its vector of ordered eigenvalues. Recall that
 	\[\lambda_*(dX)=\frac{n!\pi^{\frac{n^2}{2}}}{2^n\Gamma_n(\frac{n+2}2)}\prod _{1\leq i<j\leq n} (\lambda_i-\lambda _j) \prod_{i=1}^n d\lambda_i \]
 	Consider the integral \[I_n(s)=\int_{\Delta_n} \prod_{i<j} (\lambda_i-\lambda_j) \prod_i |\lambda_i|^sd\lambda\]
 	defined initially for $\Re(s)$ sufficiently large. Then  \[D_n^{\textbf{abs}}(s)=\frac{n!\pi^{\frac{n^2}{2}}}{2^n\Gamma_n(\frac{n+2}2)}I_n(s)\] and it is easy to see that $I_n(s)$ is a rational function of $s$: denoting $(\lambda,\mu)\in\R^a\times\R^b$, we may write
 	\[I_n(s)= \displaystyle\sum_{a+b=n}\int_{\Delta^+_a\times \Delta^+_b}  \prod_{i<j}(\lambda_i-\lambda_j)\prod_{i<j}(\mu_i-\mu_j)\prod_{i.j}(\lambda_i+\mu_j)\prod_i\lambda_i^s\prod_i \mu_i^s d\lambda d\mu \]
 	and then repeatedly apply Fubini's theorem to compute the integrals.
 	
 	Observe that for all integers $k\geq0$,
 	\[I_n(2k)= (-1)^{n\choose 2} f_n(2k+1,2k+2,\dots, 2k+n )= \frac {1}{p_n(2k)} \] where $p_n$ is the polynomial appearing in the statement. Since $I_n(s)$ is a rational function, we deduce that \[I_n(s)=\frac{1}{p_n(s)}\]
 	for all $s$, concluding the proof.
 	
 	Similarly, noting that for $s=2k+1$, $D^{\textbf{sgn}}_{n}(2k+1)$ coincides with a multiple of $f_{n}(2k+2,\dots,2k+1+n)$, we obtain the result also in that case.

 	\endproof

 	\begin{Remark}\label{rem:sign_det}
 		It is obvious directly from its definition that $D^{\textbf{\emph {sgn}}}_{2m+1}(s)=0$.
 	\end{Remark}
 	
 	Next let us record the implications of Proposition \ref{prop:mixed_signs}. Define for later use the set $\Delta_{a,n-a}=\{1\geq x_1\geq\dots\geq x_a\geq 0\geq x_{a+1}\geq \dots\geq x_{n}\geq -1\}\subset\R^n$.

 	\begin{Proposition}\label{prop:mixed_signs_applied}
 		Denote $m=\lfloor\frac{n}{2}\rfloor$. 
 		
 		\[D_n^{\textbf{\emph{cos}}}(s)+\sqrt{-1}D_n^{\textbf{\emph{sin}}}(s)=\frac{n!\pi^{\frac {n^2}{2}}}{2^{n}\Gamma_n(\frac{n+2}{2})} (-1)^{n-m\choose 2}\delta_n(m) \prod_{j=1}^n\frac{1}{s+j}\prod_{\displaystyle\scriptsize\begin{array}{c} 1\leq i<j\leq n\\ i\equiv j\mod 2 \end{array}}\frac{j-i}{2s+i+j} \]
 	\end{Proposition}
 	
 	\proof
 	
 	For $s=\frac{4k+1}{2}$ one has 
 	
 	\[\sum_{a+b=n} \sqrt {-1}^b \int_{\displaystyle\Delta_{a, b}} \prod_{j=1}^n|x_j|^{s}\prod_{i<j}(x_i-x_j)dx=(-1)^{n\choose 2}f_n(s+1,\dots,s+n).   \]
 	
 	The left hand side is a rational function of $s$, while by Proposition \ref{prop:mixed_signs} the right hand side is a complex multiple of the inverse of a real polynomial of $s$. Thus equality holds for all values of $s$, and it remains to adjust the order of the arguments in $f_n$ to arrive at the formula.
 	
 	\endproof

	\begin{Corollary}\label{cor:real_imaginary} Assume $s>0$ is a strict half-integer, and denote $R_n(s)=f_n(s+1,\dots, s+n)$. One has the following possibilities.
		
		\begin{itemize}
			\item If $n\equiv 1\mod 2$ then $\Re R_n(s)\neq0, \Im R_n(s)\neq 0$.
			\item If $n\equiv 0\mod 4$ then $\Re R_n(s)\neq 0, \Im R_n(s)= 0$.
			\item If $n\equiv 2\mod 4$ then $\Re R_n(s)=0, \Im R_n(s)\neq 0$.
		\end{itemize}
	\end{Corollary}
 	
 	\begin{Corollary}\label{cor:D_non_zero}
 		The non-zero $D^\epsilon_n(s)$, namely:  $D^{\textbf{\emph{abs}}}_n(s)$ for $n\in\mathbb N$, $D^{\textbf{\emph {sgn}}}_{2m}(s)$ for $m\in\mathbb N$, $D^{\textbf{\emph {cos}}}_n(s)$ for $n\not\equiv 2\mod 4$ and $D^{\textbf{\emph{sin}}}_n(s)$ for $n\not\equiv 0\mod 4$, are all inverse polynomials. In particular, they are non-vanishing at all $s\in\mathbb C$.
 	\end{Corollary}
 	\begin{Remark}
 		This is also the case for $D_n^+(s)$, $n\in\mathbb N$.
 	\end{Remark}
 	To prove Propositions \ref{prop:A5} and \ref{prop:mixed_signs}, we will need a few identities.
 	\begin{Lemma} 		
			For $n\geq 1$ and any $a_j\in\mathbb C$, $1\leq j\leq n$, the following identities hold.
		
			\begin{equation}\label{eq:residues}\sum_{j=1}^n a_j\prod_{\scriptsize \begin{array}{c}1\leq i\leq n\\i\neq j\end{array}}\frac{a_i+a_j}{a_i-a_j} = (-1)^{n-1}\sum_{j=1}^n a_j.\end{equation}
			
			If $n=2m$ then
			\begin{equation}\label{eq:residues1}\sum_{j=m+1}^n\frac {\displaystyle\prod_{i=1}^m (a_j+a_i)}{\displaystyle\prod_{\scriptsize \begin{array}{c} m<k\leq n\\ k\neq j \end{array}}(a_j-a_k)}=\sum_{i=1}^n a_i \end{equation}
			
			If $n=2m-1$ then
			\begin{equation}\label{eq:residues2}
			\sum_{j=1}^m  a_j \frac{\displaystyle\prod_{k=m+1}^n(a_j+a_k)}{\displaystyle\prod_{\scriptsize\begin{array}{c}1\leq i\leq m\\ i\neq j\end{array}} (a_j-a_{i})} =\sum_{j=1}^n a_j 
			\end{equation}

 	\end{Lemma}
 	\proof
 		Consider the meromorphic differential form 
 		\[\omega=\prod_{j=1}^n \frac{z+a_j}{z-a_j}dz \]
 		on $\mathbb CP^1$. It has residues \[\Res_{a_j}\omega=2a_j(-1)^{n-1}\prod_{i\neq j}\frac{a_i+a_j}{a_i-a_j}\] and 
 		\[\Res_\infty\omega= -\Res(\prod_{j=1}^n \frac{1+a_iw}{1-a_{i}w}\frac{dw}{w^2},0) =-2\sum_{j=1}^n a_j.\]
 		The sum of the residues of $\omega$ vanishes, concluding the proof of eq. (\ref{eq:residues}).
 		
 		Next consider the form $\omega= \prod_{i=1}^m \frac {z+a_i}{z-a_{m+i}}dz$. For $m+1\leq j\leq n$, \[\Res_{a_{j}}\omega=\frac {\displaystyle\prod_{1\leq i\leq m} (a_j+a_i)}{\displaystyle\prod_{\scriptsize \begin{array}{c} m<k\leq n\\ k\neq j \end{array}}(a_j-a_k)}\] while 
 		\[\displaystyle  \Res_\infty \omega = -\Res(\prod_{i=1}^m \frac{1+a_iw}{1-wa_{m+i}}\frac{dw}{w^2},0)=-\sum_{i=1}^n a_i. \]
 		Summing up the residues, we obtain eq. (\ref{eq:residues1}).
 		
 		Finally, eq. (\ref{eq:residues2}) follows immediately from eq. (\ref{eq:residues1}) by switching the lists $\{a_1,\dots, a_m\}$ and $\{a_{m+1},\dots, a_n\}$, and inserting $0$ into one of the lists.

 	\endproof
 	\begin{Remark}
		 Identity (\ref{eq:residues}) is used by DiPippo-Howe in their proof of Proposition \ref{prop:dipippo_howe}.
 	\end{Remark}

 	\noindent\textit{Proof of Proposition  \ref{prop:A5}}.\\ 
 	One easily computes that $f_1(e)=\frac{1-(-1)^e}{e}$, thus establishing the assertion for $n=1$. We then proceed by induction.
 	
 	Let us first derive a recursive relation.
 	Recall the set $\Delta_{a,n-a}=\{1\geq x_1\geq\dots\geq x_a\geq 0\geq x_{a+1}\geq \dots\geq x_n\geq -1\}$ and define the corresponding integral
 	\[f_{a,b}(e_1,\dots,e_n)=\int_{\Delta_{a,b}} \det (x_i^{e_j-1})dx \]
 	It then holds that $f_n=\sum_{a=0}^n f_{a,n-a}$.
 	Represent $\Delta_{a,b}=\Delta^+_{a,b}\cup\Delta^-_{a,b}$, where
 	$\Delta^+_{a,b}=\Delta_{a,b}\cap\{|x_1|\geq |x_n| \}$,  $\Delta^-_{a,b}=\Delta_{a,b}\cap\{|x_1|\leq |x_n| \}$, and $\Delta^+_{a,b}\cap \Delta^-_{a,b}$ has measure zero. In particular, if $b=0$ (resp. if $a=0$) then $\Delta^-_{a,0}$ (resp. $\Delta^+_{0,b}$) has measure zero.
 	
 	To evaluate $\int_{\Delta_{a,b}^+}\det(x_i^{e_j-1})dx$ we make the change of variable $x_i=t_ix_1$, $2\leq i\leq n$. We also write $t_1=1$.
 	The volume element is then $x_1^{n-1}dx_1dt_2\dots dt_n$, and \[\det( x_i^{e_j-1})=x_1^{e_1+\dots_+e_n-n}\det (t_i^{e_j-1})=x_1^{e_1+\dots_+e_n-n}\sum_{j=1}^n (-1)^{j+1}\det({t_i}^{ e_\nu-1})_{i=2\dots n}^{\nu\neq j}\] while the new domain is $\{1\geq x_1\geq 0\}\times\{1\geq t_2\geq\dots\geq t_n\geq -1\}$. 
 	
 	Therefore, 
 	\[\int_{\Delta_{a,b}^+}\det(x_i^{e_j-1})dx =\frac{1}{e_1+\dots+e_n} \sum_{j=1}^n (-1)^{j+1} f_{a-1, b}(e_1,\dots,\hat e_j,\dots, e_n) \]
 	
 	Repeating this computation for $\Delta_{a,b}^-$ where we use the change of variables $x_i=t_i |x_n|$, $1\leq i\leq n-1$, $x_n=-|x_n|$, we get the equation
 	
 	\[\int_{\Delta_{a,b}^-}\det(x_i^{e_j-1})dx =  \frac{1}{e_1+\dots+e_n}\sum_{j=1}^n (-1)^{j+n+e_j-1}f_{a, b-1}(e_1,\dots,\hat e_j,\dots, e_n) \]
 	
 	Now we sum those two integrals, and then sum over all $a+b=n$. On the right hand side, for each $j$ the multiple with $f_{\alpha,\beta}$ can correspond to either $(a,b)=(\alpha+1,\beta)$ in the $\Delta_{a,b}^+$ sum, or to $(a,b)=(\alpha, \beta+1)$ in the $\Delta_{a,b}^-$ sum. It follows that
 	
 	\begin{equation}\label{eq:recursion}f_n(e_1,\dots,e_n)=\frac{1}{e_1+\dots+e_n} \sum_{j=1}^n ( (-1)^{j+1}+(-1)^{j+n+e_j-1} )f_{n-1}(e_1,\dots,\hat e_j,\dots, e_n)\end{equation}
 	
 	Let us now check that if $N_-(e)-N_+(e)\notin\{0,1\}$ then $f_n(e)=0$.
 	
 	If $n=2m$ is even, by our assumption $|N_-(e)-N_+(e)|\geq 2$. By eq. (\ref{eq:recursion}) combined with the induction hypothesis, $f_n(e)$ will vanish unless $N_-(e)=N_+(e)+2$. Then $f_{n-1}(e_1,\dots, \hat e_j,\dots, e_n)$ can only be non-zero if $e_j$ is odd. But for such $j$, the coefficient $(-1)^{j+1}+(-1)^{j+n+e_j-1}=0$. 
 	Similarly for $n=2m-1$, by the induction hypothesis we can only have non-zero summands if $N_-(e)-N_+(e)=-1$. Then $f_{n-1}(e_1,\dots, \hat e_j,\dots, e_n)$ can only be non-zero if $e_j$ is even, and again the coefficient  $(-1)^{j+1}+(-1)^{j+n+e_j-1}=0$.
 	
 	Finally, let us assume $N_-(e)-N_+(e)\in\{0,1\}$ and establish formula (\ref{eq:super_awesome_closed_form_integral}). In the following, we write $\sum e=\sum_{j=1}^n e_j$. 
 	
 	If $n=2m$ is even, we have by the induction hypothesis
 	\[f_n(e)=\frac{2}{\sum e}\sum_{j=m+1}^{n}(-1)^{j+1}f_{n-1}(e_1,\dots,\hat e_j,\dots, e_n) \] 
 	\[=\epsilon_{n-1}\frac{2^{n}}{\sum e}\sum_{j=m+1}^n(-1)^{j+1}\frac{\displaystyle \prod_{1\leq i<i'\leq m }(e_i-e_{i'})}{\displaystyle\prod_{1\leq i\leq m }e_i}\frac{\displaystyle \prod_{\scriptsize \begin{array}{c} m<k<l\leq n\\ k,l\neq j\end{array}}(e_k-e_l)}{\displaystyle\scriptsize\prod_{\begin{array}{c} 1\leq i\leq m\\ m<k\leq n, k\neq j\end{array}} (e_i+e_k) }=\]  
 	
 	\[=\epsilon_{n-1}\frac{2^{n}}{\sum e}\sum_{j=m+1}^n(-1)^{j+1}\frac{\displaystyle \prod_{1\leq i<i'\leq m }(e_i-e_{i'})}{\displaystyle\prod_{1\leq i\leq m }e_i} \frac{\displaystyle\prod_{m<k<l\leq n}(e_k-e_l)}{\displaystyle\prod_{1\leq i\leq m<k\leq n} (e_i+e_k) }\frac {\displaystyle\prod_{1\leq i\leq m} (e_i+e_j)}{\displaystyle\prod_{\scriptsize \begin{array}{c} m<k<l\leq n\\ j\in\{k,l\} \end{array}}(e_k-e_l)} \]
 	It remains to show that
 	
 	\[\epsilon_{n-1}\sum_{j=m+1}^n(-1)^{j+1}\frac {\displaystyle\prod_{1\leq i\leq m} (e_i+e_j)}{\displaystyle\prod_{\scriptsize \begin{array}{c} m<k<l\leq n\\ j\in\{k,l\} \end{array}}(e_k-e_l)}  =\epsilon_n \sum_{j=1}^n e_j  \]
 	The left hand side is equal to 
 	
 	\[\epsilon_{n-1}\sum_{j=m+1}^n(-1)^{j+1}\frac {\displaystyle\prod_{1\leq i\leq m} (e_i+e_j)}{(-1)^{j-m-1}\displaystyle\prod_{\scriptsize \begin{array}{c} m<k\leq n\\ k\neq j \end{array}}(e_j-e_k)}\] \[=\epsilon_{n-1}(-1)^{m}\sum_{j=m+1}^n \frac {\displaystyle\prod_{1\leq i\leq m} (e_j+e_i)}{\displaystyle\prod_{\scriptsize \begin{array}{c} m<k\leq n\\ k\neq j \end{array}}(e_j-e_k)}=\epsilon_{n-1}(-1)^{m}\sum_{j=1}^n e_j\]
 	where the last equality follows from eq. (\ref{eq:residues1}). We are left to verify
 	
 	\[\epsilon_{n-1}(-1)^{m}=\epsilon_n\]
 	which readily holds.
 	
 	Now assume $n=2m-1$ is odd. Proceeding as in the even case, we are left to verify the identity
 	
 	\[ \epsilon_{n-1}\sum_{j=1}^m (-1)^{j+1}e_j \frac{\displaystyle\prod_{k=m+1}^n(e_j+e_k)}{\displaystyle\prod_{\scriptsize\begin{array}{c}1\leq i<i'\leq m\\ j\in\{i,i'\}\end{array}} (e_i-e_{i'})}=\epsilon_n\sum_{j=1}^n e_j  \]
 	or equivalently
 	\[ \epsilon_{n-1}\sum_{j=1}^m e_j \frac{\displaystyle\prod_{k=m+1}^n(e_j+e_k)}{\displaystyle\prod_{\scriptsize\begin{array}{c}1\leq i\leq m\\ i\neq j\end{array}} (e_j-e_{i})}=\epsilon_n\sum_{j=1}^n e_j  \]
    Using eq. (\ref{eq:residues2}), it remains to check that $\epsilon_{n-1}=\epsilon_n$.
 	\endproof\qed
 	
  	\noindent\textit{Proof of Proposition \ref{prop:mixed_signs}}.\\
 	For $n=0$ the equality trivially holds. We then proceed by induction. 
 	It is easy to check (carefully!) that the recursive relation (\ref{eq:recursion}) still holds. Thus
 	
 	\[f_n(e)=\]
 	\[\frac{1}{\sum e_j}\left(\displaystyle\sum_{j=1}^{m}(-1)^{j+1}(1+(-1)^{n+\frac12})f_{n-1}(\widehat e_j)+\sum_{j=m+1}^{n}(-1)^{j+1}(1-(-1)^{n+\frac12})f_{n-1}(\widehat e_j)\right) \] 
 	and by the induction hypothesis it equals
 	\[\frac{\sqrt 2}{\sum e_j\prod e_j}\cdot \]
 	
 	\[\begin{array}{c}\cdot\left(\delta_{n-1}(m-1) e^{(-1)^n\frac \pi 4\sqrt{-1}}\displaystyle\sum_{j=1}^m(-1)^{j+1}e_j\displaystyle \prod_{\scriptsize \begin{array}{c}1\leq i<i'\leq m\\ j\not\in\{i,i'\} \end{array}}\frac{e_i-e_{i'}}{e_i+e_{i'}}\prod_{m+1\leq l<l'\leq n}\frac{e_l-e_{l'}}{e_l+e_{l'}}\right.\\
 	\left.+\delta_{n-1}(m) e^{(-1)^{n+1}\frac \pi 4\sqrt{-1}}\displaystyle\sum_{j=m+1}^{n}(-1)^{j+1}e_j\prod_{1\leq i<i'\leq m}\frac{e_i-e_{i'}}{e_i+e_i'}\displaystyle \prod_{\scriptsize \begin{array}{c}m+1\leq l<l'\leq n\\ j\not\in\{l,l'\} \end{array}}\frac{e_l-e_{l'}}{e_l+e_{l'}}\right)\end{array}\]
 	which we should show is equal to 
 	\[\delta_n(m)\frac{1}{\prod e_j}\prod_{ 1\leq i<i'\leq m}\frac{e_i-e_{i'}}{e_i+e_{i'}}\prod_{ m+1\leq l<l'\leq n}\frac{e_l-e_{l'}}{e_l+e_{l'}}.\]
 	Equivalently,
 	
 	\[\begin{array}{c} \delta_{n-1}(m-1)\sqrt 2 e^{(-1)^n\frac \pi 4\sqrt{-1}}\displaystyle\sum_{j=1}^{m}(-1)^{j+1} (-1)^{m-j}e_j\displaystyle\prod_{\scriptsize\begin{array}{c}1\leq i\leq m\\ i\neq j\end{array}}\frac{e_i+e_j}{e_i-e_j} \\+\delta_{n-1}(m)\sqrt 2 e^{(-1)^{n+1}\frac \pi 4\sqrt{-1}}\displaystyle\sum_{j=m+1}^{n}(-1)^{j+1}(-1)^{n-j} e_j\displaystyle\prod_{\scriptsize\begin{array}{c}m+1\leq l\leq n\\ l\neq j\end{array}}\frac{e_l+e_j}{e_l-e_j}
 	\\
 	=\delta_{n}(m)\sum e_j\end{array}\]
 	That is, 
 	
 	\[\begin{array}{c} (-1)^{m}\delta_{n-1}(m-1)\sqrt 2 e^{(-1)^n\frac \pi 4\sqrt{-1}}\displaystyle\sum_{j=1}^{m} e_j\displaystyle\prod_{\scriptsize\begin{array}{c}1\leq i\leq m\\ i\neq j\end{array}}\frac{e_i+e_j}{e_i-e_j} \\+(-1)^{n}\delta_{n-1}(m)\sqrt 2 e^{(-1)^{n+1}\frac \pi 4\sqrt{-1}}\displaystyle\sum_{j=m+1}^{n} e_j\displaystyle\prod_{\scriptsize\begin{array}{c}m+1\leq l\leq n\\ l\neq j\end{array}}\frac{e_l+e_j}{e_l-e_j}
 	\\
 	=-\delta_{n}(m)\sum e_j\end{array}\]
 	By eq. (\ref{eq:residues}), it suffices to check the two equalities:
 	\[ \delta_{n-1}(m-1)\sqrt 2 e^{(-1)^n\frac \pi 4\sqrt{-1}}=(-1)^{m}\delta_{n-1}(m)\sqrt 2 e^{(-1)^{n+1}\frac \pi 4\sqrt{-1}}=\delta_{n}(m)\]
 	which are easy to verify for all $m,n$.
 %	\begin{itemize}
 %		\item For even $n$, \[\delta_{n-1}(m-1)\sqrt 2 e^{\frac \pi 4\sqrt {-1}} = (-1)^{m}\delta_{n-1}(m) \sqrt 2 e^{-\frac \pi 4\sqrt {-1}}=\delta_n(m)\]
 %		\item For odd $n$, 
 %		\[\delta_{n-1}(m-1)\sqrt 2 e^{-\frac \pi 4\sqrt {-1}} = (-1)^{m}\delta_{n-1}(m) \sqrt 2 e^{\frac \pi 4\sqrt {-1}}=\delta_n(m)\]
 %	\end{itemize}
 	\qed

 \section{Existence of invariant Crofton formulas}\label{sec:general_existence}
 
\subsection{Proof of Theorem \ref{thm:main_1}}
We first prove the existence of invariant Crofton formulas in a few cases.
 \begin{Theorem}\label{thm:universal_families}
 	For each of the following Crofton distributions $\mu$, $\Cr(\mu)\neq 0$.
 	\begin{enumerate}
 		\item $\mu=\mu_\textbf{\emph{abs}}(-2m)\in \mathcal M^{-\infty}(\AGr_{2m-1}(\R^{2m,2m-1}))^{tr}$.
 		\item $\mu=\mu_\textbf{\emph{sgn}}(-2m-1)\in \mathcal M^{-\infty}(\AGr_{2m}(\R^{2m+1,2m}))^{tr}$.
 		\item $\mu=\mu_{\textbf {\emph{cos}}}(-\frac{2p+1}{2})\in \mathcal M^{-\infty}(\AGr_{p-1}(\R^{p,p}))^{tr}$ when $p\not\equiv 3\mod 4$.
 		\item $\mu=\mu_{\textbf {\emph{sin}}}(-\frac{2p+1}{2})\in \mathcal M^{-\infty}(\AGr_{p-1}(\R^{p,p}))^{tr}$ when $p\not\equiv 1\mod 4$. 	
 	
 	\end{enumerate}
 \end{Theorem}
 
 \proof
 
 Consider case i). Denote $\mu=\mu_\textbf{abs}(-2m)$, $\phi=\Cr(\mu)\in \Val_{2m}^{-\infty}(\R^{2m,2m-1})$.
 We use the standard Euclidean structure to identify $\mu$ with a distribution on $\Gr_{2m-1}(\R^{2m, 2m-1})$.
 Given a function $f\in C^\infty (\Gr_{2m-1}(\R^{4m-1}))$ which is $\OO(2m)\times \OO(2m-1)$-invariant, we can write $f(E)=F(\lambda_1,\dots,\lambda_{2m-1})$ where $(\lambda_j)_{1\leq j\leq 2m-1}$ is the spectrum of $M_P(E)$. By Proposition \ref{prop:beta_distribution} one can write after a rescaling of $\mu$
 \[\langle \mu, f(E)\rangle=\int_{\Delta_{2m-1}} \prod_i \lambda_i^{s_0} \prod_{i<j}(\lambda_i-\lambda_j)\prod_i (1-\lambda_i)^{-\frac12}F(\lambda)d\lambda \]
 where $\Delta_N=\{1\geq\lambda_1\geq\dots \geq \lambda_{N}\geq -1\}$, $s_0=-2m$, and the integral is understood in the sense of analytic extension.
 
 Now let $\mathcal E_\epsilon$ be the $\OO(2m)\times \OO(2m-1)$-symmetric ellipsoid given by 
 \[\frac{1}{\epsilon^2}\sum_{i=1}^{2m}x_i^2+\sum_{i=2m+1}^{4m-1}x_i^2\leq1.\] 
 By Appendix \ref{sec:ellipsoid}, its projection to $E^P\in \Gr_{2m}(\R^{2m,2m-1})$ has volume
 \[f_\epsilon(E)=\omega_{2m}\epsilon\left(\frac{1+\epsilon^2}{2}\right)^{\frac{2m-1}{2}}\prod_{i=1}^{2m-1}\left(1-\frac{1-\epsilon^2}{1+\epsilon^2}\lambda_i(E)\right)^{\frac12}\]
 Note that \[\frac1\epsilon f_\epsilon(E)\to \omega_{2m}2^{\frac12-m}\prod (1-\lambda_i(E))^{\frac12}\] as $\epsilon\to 0$, so that $f(E):=\prod (1-\lambda_i(E))^{\frac12}$ is in the image of the cosine transform.
 Moreover, $\WF(f)= N^*S$, where $S$ are those subspaces intersecting $\R^{2m,0}$ non-generically. Indeed,
 \[ \WF(f)=M_P^*\WF(|\det(I_N-X)|^{\frac12})  \]
 while $\WF(|\det(I_N-X)|^{\frac12})$ can be deduced by a simple change of coordinates from $\WF( |\det X|^{\frac12})$, which by $\SL(N)$-invariance is the union of the conormal bundles of the strata of the set of degenerate matrices.
 
 It follows by Proposition \ref{prop:Crofton_opq_wavefront}  that 
 \[\langle \mu, f \rangle =\int_{\Delta_{2m-1}} \prod_i \lambda_i^{s_0} \prod_{i<j}(\lambda_i-\lambda_j)d\lambda \]
 By the non-vanishing of $D^{\textbf{abs}}_{2m-1}$ in Corollary \ref{cor:D_non_zero}, we conclude $\phi=\Cr(\mu)\neq 0$.
 \\\\
 The proof of case ii) is identical to i), except we use the non-vanishing of $D^{\textbf{sgn}}_{2m}$.
 
 Finally, consider cases iii) and iv). We simply apply the corresponding Crofton formula to the Euclidean ball.
 By Proposition \ref{prop:beta_distribution}, the resulting integral is a multiple of $D^{\textbf{cos}}_{p-1}(s_0)$, resp.  $D^{\textbf{sin}}_{p-1}(s_0)$.
 By Corollary \ref{cor:D_non_zero}, both are non-zero.
 
 \endproof
 
 \begin{Remark}

 By applying the $Q$-equivariant Alesker-Fourier duality $\mathbb F_Q$ (see subsection \ref{sub:fourier} for background) and invoking Lemma \ref{lem:Crofton_orthogonal_relations}, we obtain corresponding Crofton formulas for  $\Val^{-\infty}_{p-1}(\R^{p,p})^{\OO(p,p)}$ and $\Val^{-\infty}_{p-1}(\R^{p,p-1})^{\OO(p,p-1)}$.

\end{Remark}

Working towards the proof of Theorem \ref{thm:main_1}, it will be convenient to introduce the following terminology.
\begin{Definition}
 	An infinite family of triples $\mathcal F=(p_j, q_j, k_j)_{j=1}^\infty$ with $p_j,q_j,k_j\to\infty$ will be called \emph{universal} if for every $(p,q,k)$ one can find $j\geq 1$ and a collection of maps $f_\nu:\R^{a_\nu,b_\nu}\to \R^{a'_{\nu}, b'_{\nu}}$, $\nu=1,\dots, N$, such that each $f_\nu$ is either a $Q$-isometric inclusion or a $Q$-orthogonal projection, and, using $F_\nu$ to denote either $f^*_\nu$ for inclusions or $(f_\nu)_*$ for projections, the composition $F_1\circ\dots\circ F_N$ is a map between the spaces $\Val_{k_j}^{-\infty}(\R^{p_j, q_j})^{\OO(p_j, q_j)}\to \Val_{k}^{-\infty}(\R^{p, q})^{\OO(p, q)}$.  We will write $F_\mathcal F(p,q,k):\Val_{k_\mathcal F}^{-\infty}(\R^{p_\mathcal F, q_\mathcal F})^{\OO(p_\mathcal F, q_\mathcal F)}\to \Val_{k}^{-\infty}(\R^{p, q})^{\OO(p, q)}$ for any such map.
 \end{Definition}
  \noindent By Theorem \ref{thm:isomorphism}, $F_\mathcal F(p,q,k)$ is always surjective, and it is an isomorphism whenever $\min(p,q)\geq 1$.
  We will also write $F_\mathcal F(p,q,k)$ for the corresponding map between the spaces of invariant Crofton distributions, see Remark \ref{rem:wavefront_applied}.
  \begin{Lemma}\label{lem:universal}
  	If $k_j\to \infty$ and for some $0<\lambda<1$ it holds that $\lim (p_j-\lambda k_j)=\lim(q_j-(1-\lambda)k_j)=\infty$, then $(p_j, q_j, k_j)$ is universal.
  \end{Lemma}
  \proof
  
  Fix any $(p,q,k)$. Choose $j$ large enough such that the following holds: $k_j>k$, $p':=p_j-\lfloor \lambda(k_j-k)\rfloor>p$, $q':=q_j-(k_j-k-\lfloor \lambda(k_j-k)\rfloor)>q$. Then one has a projection and an inclusion
  
  \[ \R^{p_j,q_j} \xrightarrow[]{\pi}\R^{p',q'}\xleftarrow{i} \R^{p,q}\]
  yielding the chain
  \[ \Val_{k_j}^{-\infty}(\R^{p_j,q_j})^{\OO(p_j,q_j)}\xrightarrow[]{\pi_*} \Val_k^{-\infty}(\R^{p',q'})^{\OO(p',q')}\xrightarrow{i^*} \Val_k^{-\infty}(\R^{p,q})^{\OO(p,q)}   \]
  \endproof
 \noindent If follows that each of the four families in Theorem \ref{thm:universal_families} is universal.
 
 Putting together our findings, we are now ready to prove Theorem \ref{thm:main_1} which we recall.
 \begin{Theorem}
 	For all $n=p+q$ and $0\leq k\leq n$, the Crofton map
 	\[\Cr: \mathcal M^{-\infty}(\AGr_{n-k}(\R^{p,q}))^{\overline{\OO(p,q)}}\to \Val_{k}^{-\infty}(\R^{p,q})^{\OO(p,q)}\] is surjective. 
 \end{Theorem}
 \proof
 
 The Euclidean case $\min(p,q)=0$ corresponds to the classical Crofton formulas. We will assume henceforth $p\geq q\geq 1$, $1\leq k\leq n-1$.
 
 \textit{Step 1.} Let $\mathcal F$ be a universal family of triples.
 In combination with Remark \ref{rem:wavefront_applied} and Propositions \ref{prop:Crofton_restriction} and \ref{prop:Crofton_projection}, we obtain the commutative diagram 
 
 \begin{displaymath}
 \xymatrix{
 	\mathcal M^{-\infty}(\AGr_{n_\mathcal F-k_\mathcal F}(\R^{p_\mathcal F,q_\mathcal F}))^{\overline{\OO(p_\mathcal F,q_\mathcal F)}} \ar[r]^--{\Cr}\ar[d]^{F_\mathcal F(p,q,k)} &  \Val^{-\infty}_{k_\mathcal F}(\R^{p_\mathcal F,q_\mathcal F})^{\OO(p_\mathcal F,q_\mathcal F)} \ar[d]^{F_\mathcal F(p,q,k)} \\
 \mathcal M^{-\infty}(\AGr_{n-k}(\R^{p,q}))^{\overline{\OO(p,q)}}  \ar[r]^--{\Cr} &  \Val^{-\infty}_k(\R^{p,q})^{\OO(p,q)}}
 \end{displaymath}
 
 Theorem \ref{thm:isomorphism} shows the right vertical arrow is an isomorphism.
 
 Assume that for each $(p_j,q_j, k_j)\in\mathcal F$ there are $r$ invariant Crofton distributions defining linearly independent valuations. The diagram immediately implies the same holds for all $(p,q,k)$. For $r=2$, this just means the surjectivity of $\Cr$ by Theorem \ref{thm:dimension=2}. Notice that any item in Theorem \ref{thm:universal_families} provides us with a universal family with this property, with $r=1$.
 \\\\  
 \textit{Step 2.} It remains to find a universal family $\mathcal F$ such that $\Cr$ is surjective for all $(p_j,q_j,k_j)\in\mathcal F$. We will show this for $\mathcal F=\{(4j,4j,4j)\}$. Write $n=4j$.
 
	Recall that $\mathbb F_Q$ denotes the pull-back by $Q$-orthogonal complement of Crofton distributions. By step 1, one can find $\mu\in \mathcal M^{-\infty}(\AGr_{n+2}(\R^{n+2,n+2}))^{\overline{\OO(n+2,n+2)}}$ with $\phi=\Cr(\mu)\neq 0$.
	Since one of the linear combinations $\mu\pm\mathbb F_Q\mu$ is non-zero, we may assume $\mathbb F_Q\mu=\pm\mu$.

	Consider two different chains of inclusion and projections as follows:
	
	\[ \R^{n+2, n+2}\xleftarrow{i_1}\R^{n+1,n+1}\xrightarrow{\pi_1}\R^{n, n}  \]
	\[ \R^{n+2, n+2}\xleftarrow{i_2}\R^{n+2,n}\xrightarrow{\pi_2}\R^{n, n}  \]
	yielding the chains
	\[ \Val_{n+2}^{-\infty}(\R^{n+2,n+2})^{\OO(n+2,n+2)}\xrightarrow[]{i_1^*} \Val_{n+2}^{-\infty}(\R^{n+1,n+1})^{\OO(n+1,n+1)}\xrightarrow{\pi_{1*}} \Val_n^{-\infty}(\R^{n,n})^{\OO(n,n)}   \]
	\[ \Val_{n+2}^{-\infty}(\R^{n+2,n+2})^{\OO(n+2,n+2)}\xrightarrow[]{i_2^*} \Val_{n+2}^{-\infty}(\R^{n+2,n})^{\OO(n+2,n)}\xrightarrow{\pi_{2*}} \Val_n^{-\infty}(\R^{n,n})^{\OO(n,n)}   \]
	
	\begin{itemize}
		\item 	If $\mathbb F_Q \mu = -\mu$, after a rescaling one can write using Theorem \ref{thm:klain}
		\[\displaystyle \Kl(\phi)=\kappa^{\textbf{cos}}_{4j+2}=\sum_{i=0}^{2j+1}(-1)^i \kappa_{4j+2-2i, 2i}\]
		It is then easy to compute
		\[ \Kl(\pi_{1*}i_1^*\phi)=\sum_{i=1}^{2j}(-1)^i \kappa_{4j+1-2i, 2i-1} = -\kappa^{\textbf{sin}}_{4j} \]
		\[ \Kl(\pi_{2*}i_2^*\phi)=\sum_{i=0}^{2j}(-1)^i \kappa_{4j-2i, 2i} = \kappa^{\textbf{cos}}_{4j}\]
		\item 	If $\mathbb F_Q \mu = \mu$, after a rescaling $\Kl(\phi)=\kappa^{\textbf{sin}}_{4j+2}=\displaystyle\sum_{i=0}^{2j}(-1)^i \kappa_{4j+1-2i, 2i+1}$, so that \[ \Kl(\pi_{1*}i_1^*\phi)=\sum_{i=0}^{2j}(-1)^i \kappa_{4j-2i, 2i}=\kappa^{\textbf{cos}}_{4j}  \]
		\[ \Kl(\pi_{2*}i_2^*\phi)=\sum_{i=0}^{2j-1}(-1)^i \kappa_{4j-1-2i, 2i+1} = \kappa^{\textbf{sin}}_{4j} \]
		
	\end{itemize}

    It follows that $\psi_j:=\pi_{j*}i_j^*\phi$, $j=1,2$ are linearly independent. By Propositions \ref{prop:Crofton_restriction} and \ref{prop:Crofton_projection},  $\psi_j=\Cr(\pi_{j*}i_j^*\mu)$, that is given by invariant Crofton distributions. This concludes the proof.

 \endproof
 
\subsection{The centro-affine surface area} 

Denoted $\Omega_c$, it is the unique $GL(V)$-invariant upper semi-continuous valuation on convex bodies in $V=\R^p$ containing the origin in their interior, see \cite{ludwig_reitzner10}. Note that it is not translation-invariant. $\Omega_c$ can be defined as follows. For $K\subset V$ consider $K^+:=\{(x,\xi)\in K\times K^o: x\cdot \xi =1\}\subset V\times V^*$. The space $V\times V^*$ has a natural quadratic form $Q(x,\xi)=x\cdot\xi$ of signature $(p,p)$. When $K$ is smooth and has positive gaussian curvature, $K^+$ is a smooth $(p-1)$-dimensional submanifold, and the restriction of $Q$ to $K^+$ is positive-definite. We then define $\Omega_c(K)=\int_{K^+}\vol_{Q|_{K^+}}$. This definition has a natural extension to general convex bodies $K$.

Theorem \ref{thm:universal_families} allows to write rather explicit Crofton formulas for $\Omega_c(K)$ when $p\not\equiv 3\mod 4$.
\begin{Theorem}\label{thm:centro_affine}
	
	For $p\not \equiv 3\mod 4$, set 
	\[\mu_p=\left\{\begin{array}{cc}
	\mu^{p-1}_{\textbf {\emph{sin}}}(-\frac{2p+1}{2}),& p\equiv 1\mod 4\\
	\mu^{p-1}_{\textbf {\emph{cos}}}(-\frac{2p+1}{2})+(-1)^{\frac p 2}\mu^{p-1}_{\textbf {\emph{sin}}}(-\frac{2p+1}{2}),& p\equiv 0\mod 2
	\end{array}\right. \]
	The centro-affine surface area of a smooth, strictly convex body $K\subset \R^p$ is given by
	\[\Omega_c(K)=c_p\int_{\AGr_{p+1}(\R^{p,p})}\chi(K^+\cap\overline E)d\mu_p(\overline{E})\]
	for some universal constant $c_p$.
\end{Theorem}
For the sake of simplicity, we will omit the technical details justifying the applicability of the Crofton formula to $K^+$, which is non-convex.  We remark moreover that the assumptions on $K$ can be substantially relaxed.
\proof
Recall (see subsection \ref{sub:opq_classification}) the decomposition of $\OO(p,p)$-invariant valuations into $i$-even and $i$-odd.
For $p\equiv 1\mod 4$, we see by Theorem \ref{thm:universal_families} and Lemma \ref{lem:Crofton_orthogonal_relations} that $\mu=\mu^{p-1}_{\textbf {{sin}}}(-\frac{2p+1}{2})$ defines a non-trivial $i$-symmetric $(p-1)$-homogeneous valuation.
For even values of $p$, by the proof of Theorem \ref{thm:universal_families} we see that the $i$-symmetric combination $\mu_1=\mu^{p+1}_{\textbf {{cos}}}(-\frac{2p+1}{2})-(-1)^{\frac p 2}\mu^{p+1}_{\textbf {{sin}}}(-\frac{2p+1}{2})$ has $\Cr(\mu)(B)\neq 0$ for the Euclidean ball $B$, and taking the Alesker-Fourier transform we conclude that $\mu=\mu^{p-1}_{\textbf {{cos}}}(-\frac{2p+1}{2})+(-1)^{\frac p 2}\mu^{p-1}_{\textbf {{sin}}}(-\frac{2p+1}{2})$ has $\Cr(\mu)\neq 0$.
Now Theorem \ref{thm:klain} shows that for $p\not\equiv 3\mod4$, the $i$-symmetric valuation $\phi=\Cr(\mu)$ has non-trivial restriction to $X^{p-1}_{p-1,0}$. Fix $c_p\in\R$ such that $\Kl(\phi)(E)=c_p^{-1}\vol_{Q|_E}$ for $E\in X^{p-1}_{p-1,0}$.

%Observe that for a positive-definite piecewise-smooth $(p-1)$-dimensional surface $S\subset \R^{p,p}$, %$\WF(\mu_p)\cap N^*G(S)=\emptyset$ where $G(S)\subset\Gr_{p-1}(\R^{p,p})$ denotes its image under the Gauss map. 
Approximating $K^+$ by a piecewise-linear surface $S_j\to K$ and recalling that $\phi$ is Klain-Schneider continuous, we can write
\[\Omega_c(K)=c_p\lim_{j\to\infty}\phi(S_j)=c_p\lim_{j\to\infty}\int_{\AGr_{p+1}(\R^{p,p})}\chi(S_j\cap\overline E)d\mu_p(\overline{E}),\]

\[=c_p\int_{\AGr_{p+1}(\R^{p,p})}\chi(K^+\cap\overline E)d\mu_p(\overline{E}).\]

\endproof

 \section{The Crofton distribution supported on the closed orbit}\label{sec:mu_c}

 \subsection{Background from representation theory}
 
 Assume $N\leq \frac n 2$. Let us recall the description of the zonal harmonics on $\Gr_N(\R^n)=\OO(n)/\OO(N)\times \OO(n-N)$ given by James and Constantine in \cite{constantine_james}.
 
 Take $\kappa=(\kappa_1,\dots, \kappa_N)$ to be a partition of $|\kappa|:=\sum \kappa_j$ into no more than $N$ parts, $\kappa_1\geq\dots\geq \kappa_N\geq 0$. Such partitions parametrize the irreducible representations $V_\kappa$ of $\SO(2N+1)$ through their highest weight.  Moreover, the irreducible representations of $\SO(2N)$ are determined by $\kappa$ and a sign $\epsilon\in\{\pm 1\}$ through their highest weight vector $\epsilon(\kappa):=(\kappa_1,\dots,\kappa_{N-1}, \epsilon\kappa_N)$ (such vectors will be called signed partitions). They too will be denoted by $V_{\epsilon(\kappa)}$.

 The partitions $\kappa$ also parametrize the finite-dimensional irreducible polynomial representations of $\GL_\R(N)$, denoted $\rho_\kappa$. 
 
 %, for which $\rho_\kappa (\lambda I_N)=\lambda^{N|\kappa|}\cdot\Id$.
 
 The action of $g\in \GL_\R(N)$ on $X\in\Sym_N(\R)$ by $g(X)=gXg^T$ induces an action of $\GL_\R(N)$ on the polynomials on $\Sym_N(\R)$. The latter representation decomposes into a direct sum of distinct irreducible representations of the form $\rho_{2\kappa}$, see \cite{thrall} or \cite{howe}. In each $\rho_{2\kappa}$, the subgroup $\OO(N)\subset\GL_\R(N) $ has a unique one-dimensional subspace on which it acts trivially. The $\SO(N)$- and $\OO(N)$- orbits on $\Sym_N(\R)$ coincide, hence it must also be the unique subspace on which $\SO(N)$ acts trivially.
 The $|\kappa|$-homogeneous polynomial $C_\kappa(X)$ is defined to be the unique $\SO(N)$-invariant polynomial in $\rho_{2\kappa}$, normalized by $\tr(X)^{k}=\sum_{|\kappa|=k}C_\kappa (X)$. A different normalization has $C_\kappa ^*(I_N)=1$.
 
 Those polynomials satisfy many identities. We will need the following two. The first is a binomial expansion.
 
 \begin{Proposition}[Constantine \cite{constantine_hotelling}]\label{prop:plus_I}
 	\[C^*_\kappa (I_N+X)=\sum_{\sigma\leq \kappa } {\kappa \choose\sigma} C^*_\sigma (X)  \]
 \end{Proposition}
\noindent The inequality $\sigma\leq \kappa$ means $\sigma_j\leq \kappa_j$ for all $j$. The exact value of the coefficients ${\kappa \choose\sigma}$ is of no consequence for us.

The second is a generalization of the multivariate Beta integral. 
\begin{Theorem}[Constantine \cite{constantine}]\label{thm:beta_jacobi}
	\[\int_{0\leq X\leq I_N}(\det X)^s\det(I-X)^\alpha C^*_\kappa(X) dX=\frac{\Gamma_N(s+\frac{N+1}{2}, \kappa)\Gamma_N(\alpha+\frac{N+1}{2})}{\Gamma_N(s+\alpha+N+1, \kappa)} \]
	where $\Gamma_N(x)$ is given by eq. (\ref{def:gamma_n}), and $\Gamma_N(x, \kappa)=\pi^{\frac{N(N-1)}{4}}\prod_{i=0}^{N-1}\Gamma(x+\kappa_{i+1}-\frac{i}{2})$. 
\end{Theorem}

Under the action of $\OO(n)$, $C^\infty(\Gr_N(\R^n))$ decomposes into a sum of distinct irreducible representations. Their description is slightly different for even and odd values of $n$, and we will only need the odd case. Thus we assume $n=2n'+1$. Restricting to $\SO(n)$, the $\OO(n)$-irreducible components are precisely those $V_{2\kappa}$ with $2\kappa=(2\kappa_j)_{j=1}^{n'}$ for which $\kappa_{j}=0$ for $j>N$. Moreover, each $V_{2\kappa}$ has a unique one-dimensional subspace invariant under $\OO(N)\times \OO(n-N)$, spanned by the zonal harmonic $P_\kappa(X)$. Here $X=X(E)=\pi_E\pi_E^T\in\Sym_N(\R)$, where $E\in\Gr_N(\R^n)$ and $\pi_E:\R^n\to E$ is the orthogonal projection, written in the standard basis. The non-trivial spectrum of $X$ consists of $\lambda_i=\cos^2\theta_i$, the squares of the cosines of the principal angles between $E$ and the subspace $\R^N$ stabilized by $\SO(N)\times \SO(n-N)$. By $\SO(N)\times \SO(n-N)$-invariance, $\tilde P_\kappa(X)$ only depends on the spectrum $\lambda(X)$ of $X$, where $\lambda_1(X)\geq \dots\geq \lambda_N(X)$. 
 
 \begin{Theorem}[James-Constantine \cite{constantine_james}]\label{james_constantine} It holds that 
 	\[P_\kappa(X)= \sum_{\sigma\leq \kappa} \alpha_\sigma C_\sigma^* (X) \]
 	 where $\alpha_\sigma$ are given by a certain recursive relation. 
 \end{Theorem}
\noindent We will not make use of the exact values of $\alpha_\sigma$.

We will also need the description of the image of the cosine transform given by Alesker and Bernstein. Let us denote by $\Lambda^2_N(n)$ all partitions (signed partitions if $n$ is even) $\kappa=(\kappa_i)_{i=1}^N$ with $|\kappa_2|\leq 1$.
\begin{Theorem}[Alesker-Bernstein \cite{alesker_bernstein}] \label{alesker_bernstein}
	The image of the cosine transform $T_N:C^\infty(\Gr_N(\R^n))\to C^\infty(\Gr_N(\R^n))$  consists of those representations $V_{2\kappa}$ of $\SO(n)$ with $\kappa\in\Lambda^2_N(n)$.
\end{Theorem}
\noindent Note that if $\sigma\leq \kappa$ are partitions and $\kappa\in\Lambda^2_N(n)\Rightarrow \sigma\in\Lambda^2_N(n)$.

\subsection{The closed orbit}

Recall that $1\leq q\leq p$, and assume $k\leq \frac{n}{2}$.  Denoting $N=\min(k,n-k, q)$, the unique closed $\OO(Q)$-orbit in $\Gr_k(V)$ is $X^k_c:=X^k_{k-N,0}$.
By Theorem \ref{thm:1}, an $\OO(Q)$-invariant Crofton distribution supported on $X^k_{c}$ exists if and only if $n\equiv N\mod 2$, and in the latter case it is uniquely defined up to a scalar multiple by Theorem \ref{thm:1}.

For $q=1$, it follows from Theorem \ref{alesker_bernstein} that $\mu_c$ defines a non-trivial valuation. Moreover, in the next section we will see that the same is true when $q=2$. In general however it might define the zero valuation.

We consider the standard Euclidean and $(p,q)$ forms in $\R^n=\R^{p,q}$. Recall the generalized sections $f_a(s)$ which we identify with distributions through the Euclidean structure. The distribution $\mu_c$ is given by the leading Laurent coefficent of $f_{k}(s)$ around $s_0=-\frac{n+1}{2}$. Set $\alpha=\frac{|q-k|-1}{2}$, $\beta=\frac{p-k-1}{2}$.  Given  $f\in C^\infty(\Gr_{k}(\R^n))$,  by Theorem \ref{thm:1}, $\langle \mu_c,f\rangle$ is the $\lfloor\frac{N+[s_0\in \mathbb Z]}{2} \rfloor$-th Laurent coefficient at $s_0$ of

\[ \int_{1\geq \lambda_1\geq \dots\geq \lambda_n\geq 0}\prod_{i=1}^N \lambda_i^s\prod_{i<j}(\lambda_i-\lambda_j)\prod_{i=1}^N(1-\lambda_i)^{\alpha}\prod_{i=1}^N(1+\lambda_i)^{\beta}F(\lambda)d\lambda   \]
where $F(\lambda)=\int_{\OO(p)\times \OO(q)} f(gE_\lambda)dg$, and $E_\lambda\in\Gr_k(V)$ is any fixed subspace such that the  eigenvalues of $M_P(E_\lambda)$ are $\lambda=(\lambda_i)$. This is up to constant the same as

\[\int_{0\leq M\leq I} (\det M)^s\det(I-M)^\alpha \det(I+M)^\beta F(\lambda (M))dM   \]

Let us restrict to the the cases where $\beta=0$, that is $p=k+1$. Together with the condition $n\equiv N\mod 2$, it leaves 3 possibilities for $(p,q;k)$: $(2m+1,2m+1;2m)$, $(2m+1, 2m-1; 2m)$, $(2m, 2m-1;2m-1)$.  Here we only consider the last case, which has the additional property that only one Grassmannian $\Gr_{2m-1}(\R^{4m-1})$ comes into play. Moreover, $n=4m-1$ is odd, simplifying the representation-theoretic relationship between $\SO(n)$ and $\OO(n)$. We prove 

\begin{Theorem}\label{thm:mu_c}
	For $m\geq 2$ and  $(p,q,k)=(2m, 2m-1, 2m-1)$, $\mu_c$ defines the zero valuation, that is, $\mu_c$ lies in the kernel of the cosine transform.
\end{Theorem}
 
  \proof

 We have $N=k=2m-1$, $n=4m-1$, $\alpha=-\frac12$, $s_0=-2m$. For $E\in\Gr_k(\R^n)$ and $\pi_E:\R^n\to E$ the orthogonal projection, one has $X(E):=\pi_E\pi_E^T=\frac{1}{2}(I_N+M_P(E))$, see Proposition \ref{prop:beta_distribution}.
 
 For a convex body $K\subset\R^n$, let $P_K(E)$ denote the $N$-volume of its Euclidean-orthogonal projection to $E\in\Gr_N(\R^n)$.
 Set $F_K(E)=\int_{g\in \OO(N)\times \OO(n-N)} P_K(gE)dg$. By Theorem \ref{alesker_bernstein} combined with Theorem \ref{james_constantine} and Proposition \ref{prop:plus_I}, \[F_K(E)=\sum_{\kappa \in\Lambda_N^2(n)} c_\kappa P_\kappa(E)=\sum_{\kappa \in\Lambda_N^2(n)} c_\kappa' C^*_\kappa(X(E))=\sum_{\kappa \in\Lambda_N^2(n)} c_\kappa'' C^*_\kappa(M_P(E)).\]
 
By Constantine's Theorem \ref{thm:beta_jacobi}, one has
 
 \[\int_{0\leq M\leq I_N}(\det M)^s\det(I-M)^{-\frac 12} C^*_\kappa(M) dM=\frac{\Gamma_N(s+\frac{N+1}{2}, \kappa)\Gamma_N(\frac{N}{2})}{\Gamma_N(s+N+\frac12, \kappa)} \]
 For $\kappa\in\Lambda_N^2(n)$, the numerator has a pole of order either $m-1$ or $m$, while the denominator has a pole of order $m-1$. Thus we have at most a simple pole at $s_0$.
 
 Recall that $\langle \mu_c, F_K\rangle$ is the Laurent coefficient of order $\lfloor\frac{N+[s_0\in \mathbb Z]}{2} \rfloor=m$. It follows that for $m\geq 2$, $\mu_c$ defines the trivial valuation.
\endproof
 
 \section{Signature $(p,2)$}\label{sec:q=2}
 	
In this section $n=p+2$. We will work with the standard forms on $\R^n=\R^{p,2}$, and frequently identify the Crofton distributions on $\Gr_k(\R^{p,2})$ with generalized functions using the Euclidean structure, denoted $P$.

We will give an explicit basis of the $\OO(p,2)$-invariant valuation given by Crofton formulas. We will assume $2\leq k\leq p$, since for $k=1, n-1$ the cosine transform is known to be an isomorphism, and the two linearly independent Crofton distributions provided by section \ref{sec:muro} will produce a basis of valuations.

Let us list the invariant Crofton distributions we will be making use of, which are provided by Theorems \ref{thm:muro} and \ref{thm:1}, with a description of their supports. Set $s_0=-\frac{p+3}{2}$. The notation here is chosen to help keep track of the supports - it is the closure of the union of the orbits $X^k_{a,b}$ over all pairs $(a,b)$ that appear in the notation, and similarly for the closed orbit $X^k_c$. If no indices appear, the support is $\Gr_k(\R^{p,2})$. Analytic extensions are denoted by $\sigma$, and residues by $\mu$.
\begin{itemize}
	\item  If $s_0\in \mathbb Z$, we have $\mu_{k-1,0}:=\Res_{s_0}f_{k}(s)$, $\mu_{k-2,1}:=\Res_{s_0}f_{k-2}(s)$.
\item If moreover $s_0$ is even, $|\sigma|^{s_0}:=\mu_\textbf{abs}(s_0)$ is well-defined, while if $s_0$ is odd we have $\sign(\sigma)|\sigma|^{s_0}:=\mu_\textbf{sgn}(s_0)$.

\item If $s_0\notin\mathbb Z$, we have $|\sigma|_{k,0}^{s_0}-|\sigma|_{k-2,2}^{s_0}:=\mu_{\textbf {{cos}}}(s_0)$ and $|\sigma|_{k-1,1}^{s_0}:=\mu_{\textbf {{sin}}}(s_0)$, as well as $\mu_c$.
\item We will also make use of the notation $|\sigma|^s_{a,b}:=f_a(s)$ when $a+b=k$.
\end{itemize}

%Let $\mu$ be a generalized function on $\Gr_k(\R^n)$. If we can choose $\SO(n-1)\subset \SO(n)$ such that the average $\int_{\SO(n-1)}g^*\mu dg$ is non-vanishing, it will follow that $\mu$ does not lie in the kernel of the cosine transform.
	
	Let $L_k^\alpha(v)\subset \Gr_k(\R^n)$ denote the subspaces forming an angle $\alpha$ with the vector $v$. In particular, $L_k^{0}(v)=\{E: v\in E\}=\Gr_{k-1}(v^P)$, and $L_k^{{\pi}/{2}}(v)=\Gr_k(v^P)$. We fix an $\SO(n-1)=\Stab(v)$-invariant probability measure on each $L^\alpha_k(v)$
    
 	\begin{Lemma}
 		 Let $\mu$ be any $\OO(p,2)$-invariant Crofton distribution, and assume that $v^P=v^Q$. Then it holds for $L_k^{\alpha}(v)$ with both $\alpha=0,\frac{\pi}{2}$ that $\WF(\mu)\cap N^* L_k^\alpha=\emptyset$.
 	\end{Lemma}
 	\proof
	  Those are particular cases of Proposition \ref{prop:Crofton_opq_wavefront}.
	\endproof
	
	 It follows one can restrict $\mu$ as a generalized function to $L^\alpha_k$. We will then evaluate the restriction on the constant function $1$. In other words, we will be evaluating integrals of the form $\int_{L^\alpha}\mu$. They will be given by certain values of the ordinary hypergeometric function that we now describe.

	For arbitrary $a, b\in\R$, define the meromorphic in $s\in\C$ function $u(s, a ,b)=\int_{0}^{1}x^s(1+x)^a(1-x)^bdx$. This is in fact an instance of the thoroughly studied ordinary hypergeometric function, namely $u(s,a,b)={}_2F_1(-a, s+1;b+s+2;-1)$. For the sake of completeness, we state the facts that we will use. 
	
	\begin{Lemma}\label{lem:many_formulas} Let $a, b$ be arbitrary, $m\in\mathbb N$.
	\begin{equation} u(s, a, b)=u(s+1, a-1, b)+u(s, a-1, b)=\sum_{j=0}^m{m\choose j}u(s+j, a-m,b)\end{equation}
	\begin{equation} u(s, a, b)=u(s, a, b-1)-u(s+1, a, b-1)=\sum_{j=0}^m(-1)^j{m\choose j}u(s+j, a,b-m)\end{equation}
	\begin{equation} u(s, 0, 0)=\frac{1}{s+1}\end{equation}
	\begin{equation} u(s, 0, b)=B(s+1, b+1)\end{equation}
	\begin{equation} u(s, a, a)=\frac12 B(\frac{s+1}{2}, a+1)\end{equation}
	\begin{equation}\label{eq:sa0}
	u(s, a,0)=\frac{2^{a+1}}{s+1}-\left(1+\frac{a+1}{s+1}\right)u(s+1,a, 0)
	\end{equation}
	\begin{equation} \Res_{s=-m}u(s, a, b)=\frac{1}{(m-1)!}\left.\frac{d^{m-1}}{dx^{m-1}}\right|_{x=0}(1+x)^a(1-x)^b\end{equation}

	\end{Lemma}
	
	\proof
	
	All equations are straightforward verifications. Let us check eq. (\ref{eq:sa0}):
	
	\[u(s,a,0)+u(s+1,a,0)=\int_0^1 x^s({1+x})^{a+1}dx=\frac{2^{a+1}}{s+1}-\frac{a+1}{s+1}u(s+1,a,0) \]
	where the last equality is obtained through integration by parts.
	
	\endproof
	
	\begin{Remark}
		It follows from eq. (\ref{eq:sa0}) that $u(-a-2, a, 0)=-\frac{2^a}{a+1}$.
	\end{Remark}
	
	We will need later the following computation:
	
	\begin{Lemma}\label{lem:diagonal_3_computation} Assume $a,b\geq-\frac12$. Then
		\[u (-a-b-3, a, b)=2^{a+b+1}\frac{a-b}{a+1}B(-a-b-2,b+1)\]
	\end{Lemma}

	\proof
	\[ u (-a-b-3, a, b)=\int_0^1 x^{-a-b-3}(1-x)^b(1+x)^adx  \]
	
	Putting $y=\frac{x}{1+x}$ the integral becomes
	\[\int_0^{\frac12}y^{-a-b-3}(1-2y)^b(1-y)dy=\int_0^{\frac12}y^{-a-b-3}(1-2y)^b(1-y)dy=\]
	\[= 2^{a+b+2}\int_0^{1}t^{-a-b-2}(1-t)^b (1-\frac{t}{2})dt =\]
	\[=2^{a+b+2}\left(B(-a-b-2,b+1)-\frac12 B(-a-b-1, b+1)\right)\]
	It remains to verify that 
	\[B(-a-b-2,b+1)-\frac12B(-a-b-1, b+1)=B(-a-b-2,b+1)(2-\frac{a+b+2}{a+1}) \]
	\[=\frac{a-b}{a+1}B(-a-b-2,b+1)\]
	\endproof
	
	\begin{Corollary}\label{cor:diagonal_3_computation}
		If $a\geq-\frac12$ is a strict half-integer and $b\geq 0$ is an integer, it follows that $u(-a-b-3,a,b)\neq 0$.
	\end{Corollary}

	\begin{Lemma}\label{lem:value_of_restriction} Consider $L^\alpha_k(v)$ with respect to $v=e_n$, $a+b=k$.
		\[\int_{L^0_k}|\sigma|^s_{a,b}=c_{p,k} \left\{	\begin{array}{cc}0,& b=0\\ u(s, \frac{p-k}{2}, \frac{k-3}{2}),& b=1
		\\u(s, \frac{k-3}{2}, \frac{p-k}{2}),& b=2\end{array}\right.\]
		
		\[\int_{L^{\pi/2}_k}|\sigma|^s_{a,b}=c'_{p,k}\left\{	\begin{array}{cc}0,& b=2\\ u(s, \frac{p-k-1}{2}, \frac{k-2}{2}),& b=0
		\\u(s, \frac{k-2}{2}, \frac{p-k-1}{2}),& b=1\end{array}\right.\]
	\end{Lemma}
	The constants can be easily written explicitly, but will not be needed.	
	\proof
	Since all computations are very similar, we only present the case of $L^0_k$.
	If $E\in L^0_k$, it follows that $Q|_E$ is not positive definite since $Q(e_n)=-1$. This completes the case of $b=0$.
	Now parametrize $E\in L^0_k$ by $E=\Span(e_n)\oplus F$ with $F\in\Gr_{k-1}(\R^{p,1})$. Note that $|\det M_P(E)|=|\det M_P(F)|$.
	\\By Proposition \ref{prop:beta_distribution} we have:
	\begin{itemize}
		\item For $b=1$, $E\in X^k_{k-1, 1}\iff F\in X^{k-1}_{k-1,0}$, so 
		\[\int_{L^0_k}|\sigma|^s_{k-1,1}=c_{p,k}\int_0^1x^s(1-x)^{\frac{k-3}{2}}(1+x)^{\frac{p-k}{2}}dx \]
		\item For $b=2$, $E\in X^k_{k-2, 2}\iff F\in X^{k-1}_{k-2,1}$, so
		 \[\int_{L^0_k}|\sigma|^s_{k-2,2}=c'_{p,k}\int_{-1}^0|x|^s(1-x)^{\frac{k-3}{2}}(1+x)^{\frac{p-k}{2}}dx \]
	\end{itemize} 
	
	\endproof
		
	\begin{Lemma}\label{lem:odd_p_residue}
		Let $a,b\geq-\frac12$ be strict half-integers, and $a\neq b$. It then holds that $\Res_{-a-b-3} u(s,a, b)\neq0$.
	\end{Lemma}
	\proof
	 Denote $m=a+b+3$. By Lemma \ref{lem:many_formulas} we have 
	\[\Res_{-m} u(s,a,b)= \frac{1}{(m-1)!} ((1+x)^a(1-x)^b)^{(m-1)}(0) \]
	We are to show the latter expression is non-zero. Replacing $x$ with $-x$ if necessary, we may assume $a>b$. Rewrite \[(1+x)^a(1-x)^b=(1+x)^{a-b}(1-x^2)^b= \sum_{i,j\geq0} (-1)^j{a-b\choose i}{b\choose j}x^{i+2j}\]
	Thus we are left to show that
	\[\sum_{i+2j=m-1}(-1)^j{a-b\choose i}{b\choose j}\neq0 \]
	
	Since $a-b\in\mathbb N$, the only values of $i$ entering the sum are $0\leq i\leq a-b$, so that $j\geq \frac{1}{2}(a+b+2-(a-b))=b+1$, so in fact $j>b+1$ as $b$ is not an integer.
	It follows that the signs of ${b\choose j}=\frac{b(b-1)\dots(b-j+1)}{j!}$ alternate as $j>b+1$ increases, implying all the summands have the same sign. This concludes the proof.   
		
	\endproof

	\begin{Lemma}\label{lem:odd_p_total}
			Assume $a,b\geq0$ are half-integers (possibly integers), and $|b- a|=m\in\mathbb Z$.  
			Then \[\lim_{s\to -a-b-3}u(s, a,b)+(-1)^{m+1}u(s, b, a)\neq0\]
	\end{Lemma}
	
	\proof 
	We may assume $b=a+m$.	Put $s_0=-a-b-3=-2a-m-3$. By Lemma \ref{lem:many_formulas}, we have 
	\[u(s, a+m, a)+(-1)^{m+1}u(s, a,a+m)=\sum_{j=0}^m{m\choose j} (1+(-1)^{j+m+1})u(s+j, a, a)\]
	\[=\frac12 \sum_{j=0}^{m}{m\choose j}(1+(-1)^{j+m+1})B(\frac{s+j+1}{2}, a+1) \]
	Write 
	\[B(\frac{s+j+1}{2}, a+1) =\frac{\Gamma(a+1)\Gamma(\frac{s+j+1}{2})}{\Gamma(\frac{s+j+3}{2}+a)}. \]
	The non-zero summands are those with $j\not\equiv m\mod 2$. Now since $\frac{s_0+j+3}{2}+a=\frac{j-m}{2}$ is a strict half-integer, and $\frac{s_0+j+3}{2}+a\leq\frac{s_0+m+3}{2}+a=0$, it follows that the signs of both the numerator and the denominator alternate between $j$ and $j+2$. It follows all summands have the same sign, completing the proof.   
	
	\endproof
	
	We will need the classification of invariant Crofton distributions in $\R^{p,2}$ with $p$ even.
	
	\begin{Proposition}\label{prop:crofton_classification}
			Let $p$ be even, denote $s_0=-\frac{p+3}{2}$. The $\OO(p,2)$-invariant Crofton distributions are spanned by the basis $\{|\sigma|^{s_0}_{k,0}-|\sigma|^{s_0}_{k-2,2}, |\sigma|^{s_0}_{k-1,1}, \mu_c\}$.
		\end{Proposition}
	\proof
	By Proposition \ref{prop:generalCroftonBound}, there is at most a 4-dimensional space of invariant Crofton distributions. Assume it is in fact 4-dimensional. Let $\mu$ be an invariant distribution independent of the listed three. By subtracting certain multiples of $|\sigma|^{s_0}_{k,0}-|\sigma|^{s_0}_{k-2,2}$, $|\sigma|^{s_0}_{1,1}$, Proposition \ref{prop:generalCroftonBound} allows us to assume $\supp \mu= \overline {X^k_{k,0}}$. Fix a decomposition $\R^{p,2}=\R^{p-2, 0}\oplus \R^{2,2}$ and a subspace $E_0\in\Gr_{k-2}\R^{p-2,0}$. Identify $\Gr_2{\R^{2,2}}$ with a submanifold $X\subset \Gr_k{\R^{p,2}}$ through the embedding $i(F)=F\oplus E_0$. 
	We consider $\mu$ as a generalized function on $\Gr_k\R^{p,2}$. Applying Proposition \ref{prop:Crofton_opq_wavefront}, we may consider  $\nu=i^*\mu\in C^{-\infty }(\Gr_2{\R^{2,2}})$. It holds that $\supp\nu=\overline {X^2_{2,0}}$, and  $g^*\nu=\psi_g(F)^{s_0} \nu$ for all $g\in \OO(2,2)$.
	
	Consider the meromorphic in $s$ generalized function $f_2(s)\in C^{-\infty}(\Gr_2\R^{2,2})$ constructed in Proposition \ref{prop:muro_pullback} (identified with functions), which satisfies $g^*f_2(s)=\psi_g(F)^sf_2(s)$. By Theorem \ref{thm:muro}, $f_2$ has Laurent series given by 
	\[f_2(s)= \frac{h_{-1}}{s-s_0}+h_0+\dots  \]
	where $h_{-1}$ is supported on $X^2_c$ and satisfies $g^* h_{-1}=\psi_g(F)^{s_0}h_{-1}$, and $\supp h_0=\overline{X^2_{2,0}}$.
	Moreover, it follows from comparing the Laurent series of the equation $g^*f_2(s)=\psi_g(F)^sf_2(s)$  that \[g^*h_0=\psi_g(F)^{s_0}h_0+\psi_g(F)^{s_0}\log\psi_g(F)h_{-1}.\]
	In particular, the restriction of $h_0$ to the complement of $X^2_c$ satisfies $g^*h_0=\psi_g(F)^{s_0}h_0$. It follows from Proposition \ref{prop:generalCroftonBound} that after an appropriate rescaling of $\nu$, $\delta=h_0-\nu$ has $\supp (\delta)=X^2_c$.  
	Moreover, $\delta$ satisfies 
	\begin{equation}\label{eqn:quasi_homogeneous}g^*\delta -\psi_g(F)^{s_0}\delta=\psi_g(F)^{s_0}\log\psi_g(F)h_{-1} \end{equation}
	In particular, $\delta$ is $\OO(2)\times \OO(2)$-invariant as $\psi_g(F)=1$ for $g\in \OO(2)\times \OO(2)$.
	Now fix $F_0\in X^2_c$, recall that $N_{F_0}X^2_c=\Sym^2F_0^*$, and by Witt's extension theorem $\Stab(F_0)$ acts by $\GL(F_0)$ on $F_0$. By $\OO(2)\times \OO(2)$-invariance and its transitive action on $X^2_c$, the wavefronts of $\delta$ and $h_{-1}$ lies within $N^*_{F_0}{X^2_c}$. We may therefore restrict to any 3-dimensional submanifold through $F_0$ which is transversal to $X^2_c$. Choose such a submanifold $S$ near $F_0$ and $g_\lambda\in \Stab(F_0)$ stabilizing $S$ locally such that $g_\lambda|_{N_{F_0}X^2_c}=\lambda\cdot Id$. In part icular, $\psi_{g_\lambda}(F_0)=|\det {g_\lambda}|_{F_0}|^{-2}=\lambda^{-4}$.
	Choose a coordinate chart in $\Gr_2(\R^{2,2})$ near $F_0$ and identify $N_{F_0}X^2_c=T_{F_0}S=\R^3$.
	Note that for $f\in C^\infty_c(\R^3)$, $\langle g_\lambda^*\delta, f\rangle $ and $\langle \psi_{g_\lambda}(F_0)^{s_0} \delta,f\rangle$ are rational functions of $\lambda$. The right hand side of equation (\ref{eqn:quasi_homogeneous}) is however transcendental in $\lambda$, a contradiction. Thus $\mu$ does not exist.

	\endproof
	
	\begin{Remark} It is likely one can extend this method to show that the Crofton distributions constructed in Proposition \ref{thm:1} span the space of invariant distributions.
		
		We remark also that $h_0$ is the pull-back to the Grassmannian of what is sometimes called a quasi-homogeneous generalized function on $\Sym_2(\R)$.
	\end{Remark}
	 We are now ready to prove Theorem \ref{thm:main_q=2}.
	\begin{Theorem}\label{thm:q=2}
		For all values of $p\geq 2$, $2\leq k\leq p$, the $\OO(p,2)$-invariant, $(p+2-k)$-homogeneous valuations are given by the following Crofton distributions:
		\[\left\{\begin{array}{lllll} \Span\{|\sigma|^{s_0}_{k,0}-|\sigma|^{s_0}_{k-2,2}, \mu_c\},& p\equiv k\equiv0\mod 2
		\\  \Span\{|\sigma|^{s_0}_{k-1,1}, \mu_c\},& p\equiv 0, k\equiv 1\mod2
		\\ \Span\{|\sigma|^{s_0}, \mu_{k-2+\epsilon,1-\epsilon}\},&  p\equiv 1 \mod 4, k\equiv\epsilon\mod 2
		\\\Span\{\sign\sigma|\sigma|^{s_0}, \mu_{k-2+\epsilon,1-\epsilon}\},&  p\equiv 3 \mod 4, k\equiv\epsilon\mod 2
			
		\end{array}\right.\]
	\end{Theorem}
	\proof
	
	Note that by Lemma \ref{lem:Crofton_orthogonal_relations}, each row is invariant under the $\OO(p,q)$-equivariant Alesker-Fourier duality (we refer to subsection \ref{sub:fourier} for background), except that the value of $k$ gets replaced by $p+2-k$. Thus, given a value of $k$ at any particular case, we only need to verify the statement for either of the values $(k, p+2-k)$. 
	
	Fix the standard $\SO(n-1)\subset \SO(n)$ fixing $e_n$. The $\SO(n)$-orbit of an $\SO(n-1)$-invariant element in $C^{-\infty}(\Gr_k(\R^n))$ defines an irreducible representation of $\SO(n)$ by spherical harmonics, which correspond to a partition $(\kappa)$ of length 1. By Theorem \ref{alesker_bernstein}, the kernel of the cosine transform $T_k$ intersects the spherical harmonics trivially. It follows that if for some generalized functions $(\mu_j)_{j=1}^r$ the averages $[\mu_j]_{\SO(n-1)}:=\int_{\SO(n-1)}g^*\mu_j dg$ are linearly independent, then so are $(T_k\mu_j)_{j=1}^r$.
	
	We will write $L_k^\alpha=L_k^\alpha(e_n)$, and repeatedly use Lemma \ref{lem:value_of_restriction} without mention. Note that when  the restriction of $\mu$ to ${L_k^\alpha}$ is well-defined then $\int_{L_k^\alpha}\mu=\int_{L_k^\alpha}[\mu]_{\SO(n-1)}$, for $\alpha=0,\frac \pi 2$, and those two linear functionals of $\mu$ are linearly independent.
	
	%In turn, if $\mu$ is an $O(p,2)$-invariant
		
%	by computing the integral of their restriction to $L_k^\alpha$ for $\alpha=0,\frac \pi 2$.
%	We thus will be able to conclude that the non-averaged Crofton distributions define linearly independent valuations.

	\noindent\textit{Step 1. The case of $p\equiv 3\mod 4$, or $p\equiv 1\mod 4$ and $k\neq \frac{p+1}{2},\frac{p+3}{2}$.}\\
     Let $\mu_1,\mu_2$ be the corresponding Crofton distributions appearing in the third or fourth row.
     By the previous paragraph, it suffices to show that $\left(\int_{L_k^\alpha}\mu_j\right)_{\alpha=0,\frac \pi 2}^{j=1,2}$ is an invertible matrix. This in turn is a consequence of the following observations.
	
	\begin{itemize}
		\item It holds that $\supp\mu_{k-1,0}\cap L_k^0=\emptyset$ as well as $\supp \mu_{k-2,1}\cap L^{\pi/2}_k=\emptyset$.
		
		\item If $k\neq \frac{p+3}{4}$ is even, it follows from Lemma  \ref{lem:odd_p_residue} that $\int_{L^0_k}\mu_{k-2,1}\neq0$.
		\item If $k\neq \frac{p+1}{2}$ is odd, Lemma  \ref{lem:odd_p_residue} implies  $\int_{L^{\pi/2}_k}\mu_{k-1,0}\neq0$.
		\item 	For $p\equiv 1\mod 4$ and $k$ even, it follows from Lemma \ref{lem:odd_p_total} with $a=\frac{p-k-1}{2}$, $b=\frac{k-2}{2}$ that $\int_{L^{\pi/2}_k} |\sigma|^{s_0}\neq 0$.
		\item 	For $p\equiv 1\mod 4$ and $k$ odd, Lemma \ref{lem:odd_p_total} with $a=\frac{p-k}{2}$, $b=\frac{k-3}{2}$ implies $\int_{L^{0}_k} |\sigma|^{s_0}\neq 0$.
		\item  For $p\equiv 3\mod 4$ and $k$ even, by Lemma \ref{lem:odd_p_total}, $\int_{L^{\pi/2}_k}  \sign\sigma|\sigma|^{s_0}\neq 0$.
		\item For $p\equiv 3\mod 4$ and $k$ odd, by Lemma \ref{lem:odd_p_total}, $\int_{L^{0}_k} \sign\sigma|\sigma|^{s_0}\neq 0$.
	\end{itemize}
	\textit{Step 2. Even $p$ - part 1.}  We will show that $\mu_c$ defines a non-trivial valuation.
	\\Assume first that either $ p\not\equiv 0\mod 4$ or $k\neq \frac{p+2}{2}$.
	Fix an isometric embedding $j:\R^{p,2}\to \R^{p+1,2}$. Consider the Crofton distribution $\mu^{k+1}\in\mathcal M^{-\infty}(\AGr_{k+1}\R^{p+1,2})^{tr}$ given by $\mu^{k+1}:=\mu_{k-\epsilon,\epsilon}$, where $k\equiv \epsilon\mod 2$.
	We proved in step 1 that $\Cr(\mu)\neq 0$. 
	
	Assume for clarity $\epsilon=0$. By Muro's theorem \ref{thm:muro}, $\supp \mu^{k+1}\subset X^{k+1}_{k,0}\cup X^{k+1}_{k-1,0}$. Let $j^*$ denote the restriction of Crofton distributions, which is essentially the push-forward under the intersection with $\R^{p,2}$ map.
	Observe  that the intersection of $\R^{p,2}$ with a non-negative semi-definite subspace is again non-negative semi-definite, so that \[\supp j^*\mu^{k+1}\subset \overline {X^{k}_{k,0}}=X^{k}_{k,0}\cup X^{k}_{k-1,0}\cup X^{k}_{k-2,0}.\]
	
	Now by Proposition \ref{prop:crofton_classification}, we conclude $j^*\mu^{k+1}=c\mu_c$ for some constant $c$.
	Since $0\neq j^*\Cr(\mu^{k+1})=\Cr(j^*\mu^{k+1})=c \Cr(\mu_c)$, we conclude $ \Cr(\mu_c)\neq 0$.
	
	Finally, in the event that $p\equiv 0\mod 4$ and $k= \frac{p+2}{2}$, we embed $j:\R^{p,2}\to \R^{p+3,2}$. 
	We have the Crofton distribution $\mu_{k+1,1}\in\mathcal M^{-\infty}(\AGr_{k+3}\R^{p+3,2})^{tr}$, which by step 1 defines a non-trivial valuation since $k+3$ is even, and it is supported on the non-positive subspaces. We now proceed as in the previous case.
	\\\\
	\noindent \textit{Step 3. Even $p$ - part 2.} Observe that the restriction of $\mu_c$ to $L^0_k$ vanishes, as the support $X^k_c$ of $\mu_c$ is disjoint from $L^0_k$. Denoting the other Crofton distribution (in the corresponding row of the statement) by $\mu$, we will show that $\int_{L^0_k}\mu\neq 0$. Then $T_k\mu_c, T_k\mu\neq 0$, and for $a,b\neq0$, $\int_{L_k^0}(a\mu+b\mu_c)=a\int_{L^0_k}\mu\neq0\Rightarrow T_k(a\mu+b\mu_c)\neq 0$. Applying Lemma \ref{lem:value_of_restriction} separately for even and odd $k$ and in each case invoking Corollary \ref{cor:diagonal_3_computation}, we establish the statement of the theorem for even $p$ and arbitrary $k$. 
	\\\\
	\textit{Step 4. The remaining case of $p\equiv 1\mod 4$, $k=\frac{p+3-2\epsilon}{2}$, $\epsilon=0,1$.}\\
    By replacing $k$ with $n-k$ we may assume $k=\frac{p+1}{2}$ is odd. By step 1, $\int_{L^0_k}|\sigma|^{s_0}\neq 0$ while $\int_{L_k^0} \mu_{k-1,0}=0$.
	
	Arguing as in step 3, it only remains to show that $\mu_{k-1,0}$ defines a non-trivial valuation. 
	Consider an inclusion $j:\R^{p,2}\to \R^{p+1,2}$. By step 2, $\mu_c$ defines a non-trivial valuation, and the restriction $j^*\mu_c$ is supported on $X^k_{k-1,0}\cup X^k_{k-2,0} $. By Proposition  \ref{prop:generalCroftonBound}, $j^*\mu_c$ coincides with a multiple of $\mu_{k-1,0}$, and by Proposition \ref{prop:Crofton_restriction} and Theorem \ref{thm:isomorphism} we conclude $\Cr(\mu_{k-1,0})\neq 0$. 
	
   \endproof

 \begin{appendix}

 \section{Projections of an $\OO(p)\times \OO(q)$-invariant ellipsoid}\label{sec:ellipsoid}
 
 Let $E\in\Gr_k(\R^n)$ form $N=\min(k,n-k,p,n-p)$ principal angles $\theta_1,\dots,\theta_{N}$ with $\R^{p}\subset\R^n$. Let $\pi_E:\R^n \to E$ be the orthogonal projection.
 Write $\lambda_i=\cos2\theta_i$, $A=\frac{a^2+b^2}{2}$, $B=\frac{a^2-b^2}{2}$. Write $\omega_k$ for the volume of the Euclidean unit ball in $\R^k$.  Set $q=n-p$,  and let $\R^p\oplus \R^q=\R^n$ be the decomposition into coordinate subspaces. Write $N_q=\max(k-q, 0)$, $N_p=\max(k-p, 0)$.
 \begin{Lemma}
  Consider the ellipsoid $\mathcal E$ given by the equation \[\frac{1}{a^2}\sum_{j=1}^px_j^2+\frac{1}{b^2}\sum_{j=1}^q x_{p+j}^2=1.\]
  Then 
 \[\vol_k(\pi_E(\mathcal E))=\omega_k a^{N_q}b^{N_p} \prod_{j=1}^N(a^2\cos^2\theta_j+b^2\sin^2\theta_j)^{1/2}=\omega_k a^{N_q}b^{N_p} \prod_{j=1}^N(A+B\lambda_j)^{1/2}\]
 \end{Lemma}
 \proof
  
 Let us denote $y=(y_j)_{j=1}^p\in\R^p$, $z=(z_j)_{j=p+1}^{p+q}\in \R^q$. Let $e_j$ be the standard basis of $\R^p$ and $f_j$ that of $\R^q$. We may choose \[E=\Span(\cos\theta_j e_j +\sin\theta_j f_j)_{j=1}^N\oplus \Span(e_{j})_{j=N+1}^{N+N_q}\oplus \Span(f_{j})_{j=N+1}^{N+N_p}\] and those vectors constitute an orthonormal basis of $E$. We have \[E^\perp=\Span\{\sin \theta_j e_j-\cos\theta_j f_j\}_{j=1}^N\oplus \Span(e_{j})_{j=N+N_q+1}^{p}\oplus \Span(f_{j})_{j=N+N_p+1}^{q}\]
  
 For $x=y+z\in\mathcal E$, $ T_x\mathcal E=\{\frac {1}{a^2}ydy+\frac{1}{b^2}zdz=0\}=(\frac {y}{a^2}+\frac{z}{b^2})^\perp$. Now $x$ is in the shadow boundary of the projection of $\mathcal E$ to $E$ if and only if $ (T_x\mathcal E)^\perp\subset E \iff \frac {y}{a^2}+\frac{z}{b^2}\perp E^{\perp}$. This amounts to the equations $z_j=\frac{b^2}{a^2}\tan\theta_j y_j$ for $1\leq j\leq N$, and $y_j=0$ for $ N+N_q<j\leq p$ and $z_j=0$ for $N+N_p<j\leq q$.
 
 Since $x\in\mathcal E$, we get the equation \[\sum_{j=1}^N \left(\frac 1 {a^2}+\frac{1}{b^2}\frac{b^4}{a^4}\tan^2\theta_j\right)y_j^2 +\frac{1}{a^2}\sum_{j=N+1}^{N+N_q}y_j^2+\frac{1}{b^2}\sum_{j=N+1}^{N+N_p}z_j^2= 1\]
 Write \[\pi_E(x)=\sum_{j=1}^N \alpha_j(\cos\theta_j e_j+\sin\theta_j f_j)+ \sum_{j=N+1}^{N+N_q}\beta_j e_j + \sum_{j=N+1}^{N+N_p}\gamma_j f_j.\]

 From $(x-\pi_E(x))\perp E$ we get:
 \begin{itemize}
 	\item  For $1\leq j\leq N$, 
 	\[(y_j-\alpha_j\cos\theta_j)\cos\theta_j +(z_j-\alpha_j\sin\theta_j)\sin\theta_j =0\]
 	so that \[y_j=\left(\cos\theta_j + \frac{b^2}{a^2} \frac{\sin^2\theta_j}{\cos\theta_j}\right)^{-1}\alpha_j=\frac{a^2}{a^2\cos\theta_j +b^2\sin^2\theta_j }\alpha_j\].
 	\item For $N+1\leq j\leq N+N_q$, $y_j=\beta_j$.
 	\item For $N+1\leq j\leq N+N_p$, $z_j=\gamma_j$.
 \end{itemize} 
 
In the orthonormal coordinates $\alpha_j, \beta_j, \gamma_j$ on $E$, the boundary of the projection of $\mathcal E$ is then given by 
\[ \sum_{j=1}^N \frac{1}{a^2\cos^2\theta_j+b^2\sin^2\theta_j}\alpha_j^2+ \frac{1}{a^2}\sum_{j=N+1}^{N+N_q}\beta_j^2+\frac{1}{b^2}\sum_{j=N+1}^{N+N_p}\gamma_j^2=1  \]

Thus the volume of the projection is 
\[\vol_k(\pi_E(\mathcal E))= \omega_ka^{N_q}b^{N_p}\prod_{j=1}^N (a^2\cos^2\theta_j+b^2\sin^2\theta_j)^{\frac12}\qedhere \]

  %\section{Pull-backs and push-forwards of Crofton distributions. by Dmitry Faifman and Thomas Wannerer}  

  \section{
  	Pull-backs and push-forwards of Crofton distributions. 
  	\\\normalfont{by Dmitry Faifman and Thomas Wannerer}
  }\label{sec:functorial}
  
  The operations of pull-back and push-forward of translation-invariant valuations under linear maps were defined by Alesker in \cite{alesker_fourier} for continuous valuations, and extended in \cite{bernig_faifman_opq} to the class of Klain-Schneider continuous valuations. Moreover, since pull-back by injection and push-forward by surjection preserve the class of smooth valuations, one obtains by Alesker-Poincar\'{e}  duality the operations of push-forward by injection and pull-back by surjection between the corresponding spaces of generalized valuations, see e.g. \cite{alesker_fourier} for definitions and details.  
  
  Since $\Cr: \mathcal M^{-\infty}(\AGr_{n-k}(V))^{tr}\to \Val_k^{+,-\infty}(V)$ is surjective (but generally not injective), one may ask if the restriction of any of those operations to even valuations can be carried out already on the level of Crofton distributions. In all cases, the answer turns out to be positive, with the caveat that in  two cases - that of pull-back by injection and push-forward by surjection - it only holds under a certain extra assumption on the wavefronts. 
  
  In the remainder of the section, $j:U\to V$ is an inclusion of linear spaces, $\pi:V\to W$ a surjection of linear spaces, $\dim V=n$, $\dim U=\dim W=m=n-d$. Subspaces of $V$ will be denoted by $E$, subspaces of $U$ and $W$ will be denoted by $F$. Note that the push-forwards actually operate on spaces of valuations twisted by dual densities. Elements of $\Val(V)\otimes \Dens^*(V)$ will be referred to as dual valuations, and similarly for Crofton distributions.

  Let us recall some basic facts concerning the functorial properties of distributions on manifolds. For details and the basics of microlocal analysis, we refer to \cite{duistermaat},\cite{guillemin_sternberg}.  
  
  Let $f:X\to Y$ be a proper map between smooth manifolds. Then there is an induced push-forward map 
  \[f_*:\mathcal M^{-\infty}(X)\to \mathcal M^{-\infty}(Y)\]
  Moreover, there is a pull-back map
  \[f^*: \mathcal M^\infty(Y)\to C^\infty(X, f^*|\omega_Y|) \]
  which extends to spaces of distributions with controlled wavefront set, namely
  \[f^*: \mathcal M_\Gamma^{-\infty}(Y)\to C^{-\infty}(X, f^*|\omega_Y|) \]
  where $\Gamma\subset T^*Y\setminus 0$ is a closed cone such that $\Gamma\cap \Ker(df^*)=\emptyset$, and $\mathcal M^{-\infty}_\Gamma(Y)=\{\mu\in\mathcal M^{-\infty}(Y): \WF(\mu)\subset\Gamma\}$.

  \subsection{Pull-back by surjection.}
  This and the next case are the simpler of the four. Given $\pi: V\to W$, there is a continuous pull-back of generalized valuations: $\pi^*:\Val_k^{-\infty}(W)\to \Val_k^{-\infty}(V)$. 
  
  Consider the proper smooth map $\pi^{-1}: \AGr_{m-k}(W)\to\AGr_{n-k}(V)$ given by $\pi^{-1}(\overline F+x)=\pi^{-1}(\overline F)+\tilde x$ for any $\overline F\in\AGr_{m-k}(W)$ and $x\in W$, $\tilde x\in V$ such that $\pi(\tilde x)=x$.
  We get a sequentially continuous push-forward map \[(\pi^{-1})_*: \mathcal M^{-\infty}(\AGr_{m-k}(W))^{tr}\to \mathcal M^{-\infty}(\AGr_{n-k}(V))^{tr} \]
  Note that under this map, Borel measures are mapped to Borel measures, while smooth measures do not in general remain smooth.
  To retain compatibility with the notation for valuations, we will write $\pi^*:=(\pi^{-1})_*$.
  \begin{Proposition}\label{prop:Crofton_pushforward_surjection}
  	Let $\mu$ be any Crofton distribution for $\phi\in \Val_{k}^{+,-\infty}(W)$. Then 
  	$\pi^*\phi=\Cr(\pi^*\mu)\in \Val_k^{+,-\infty}(V)$.
  	
  \end{Proposition}
  
  \proof
  Assume $\phi=\Cr(\mu)$ in $W$, and take $K\in \K_s(V)$. Then
  \[\pi^*\phi (K)=\phi(\pi( K))=\int_{\AGr_{m-k}(W)}\chi(\pi(K)\cap \overline F)d\mu(\overline F) \]
  \[=\int_{\AGr_{m-k}(W)}\chi(K\cap \pi^{-1}\overline F)d\mu(\overline F)=\int_{\AGr_{n-k}(V)}\chi(K\cap \overline E)d(\pi^{-1}_*\mu)(\overline E)=\Cr(\pi^{*}\mu)(K)\]
  
  \endproof
  
  \subsection{Push-forward by inclusion.}
  Given an inclusion $j:U\to V$, we get the map:
  $j_*:\Val_k^{-\infty}(U)\otimes\Dens^*(U)\to \Val_{k+d}^{-\infty}(V)\otimes \Dens^*(V)$. 
  Now $j$ extends to give a smooth map $j: \Gr_{m-k}(U)\to \Gr_{m-k}(V)$. 
  We consider the corresponding twisted Crofton distributions as generalized measures over the linear Grassmannian with values in a certain line bundle:
  \[\mathcal M^{-\infty}(\Gr_{m-k}(V), \Dens(V/E))\otimes \Dens^*(V) = \mathcal M^{-\infty}(\Gr_{m-k}(V), \Dens^*(E)) \]
  and similarly for $U$.
  
  We get a sequentially continuous push-forward of dual Crofton distributions 
  \[j_*:\mathcal M^{-\infty}(\Gr_{m-k}(U), \Dens^*(F)) \to \mathcal M^{-\infty}(\Gr_{m-k}(V),  \Dens^*(E))\]
  which respects Borel but not smooth   measures.
  
  \begin{Proposition}\label{prop:Crofton_pushforward_inclusion}
  	If $\phi\in \Val_{k}^{+,-\infty}(U)\otimes \Dens^*(U)$ and $\phi=\Cr(\mu)$, then 
  	$j_*\phi=\Cr(j_* \mu)\in \Val_{k+d}^{+,-\infty}(V)\otimes\Dens^*(V)$.
  	
  \end{Proposition}
  
  \proof
  Fix a Euclidean structure on $V$ to fix Lebesgue measures on subspaces and quotients.
  We have for $K\in\K_s(V)$
  \[j_*{\phi}(K)=\int_{V/U} \phi((K+x)\cap U)dx=\int_{V/U}dx\int_{\AGr_{m-k}(U)}\chi ((K+x)\cap \overline F )d \mu(\overline F) \]
  which is immediately seen to coincide with \[\Cr(j_*\mu)(K)=\int_{\AGr_{m-k}(V)} \chi(K\cap \overline E) d(j_*\mu)(\overline E)=\int_{\Gr_{m-k}(U)}\int_{V/F}\chi(K\cap (F+x))d\mu\]
  since $\Dens(V/U)\otimes \Dens(U/F)=\Dens(V/F)$.
  \endproof
  
  \subsection{Pull-back by inclusion.}
  We identify $U$ with $jU\subset V$. Thus assume $U\subset V$, $\dim U=m$, $\dim V=n$, $n=m+d$, and assume $k\geq d$. Consider for $E\in\Gr_{n-k}(V)$ the subspace $J(E)=E\cap U$. This gives a smooth map $J :\Gr_{n-k}(V)\setminus S\to \Gr_{m-k}(U)$, where $S$ consists of those subspaces that intersect $U$ non-generically. $S$ is stratified by locally closed submanifolds, given by $S_r=\{E:\dim(E\cap U)=n-k-d+r\}$, $1\leq r\leq d$ and $\overline{S_1}\supset\dots\supset \overline{S_d}=S_d$.
  Let $N^*S$ denote the union of the conormal bundles of the strata in $T^*\Gr_{n-k}(V)$.
  
  Let $\Gamma\subset T^*\Gr_{n-k}(V)\setminus 0$ be a closed cone. Let $\Val_{k}^{+,KS}(V;\Gamma)$ denote the class of even Klain-Schneider continuous valuations admitting a Crofton distribution with wavefront belonging to $\Gamma$. 
  
  \begin{Proposition}\label{prop:intersection_wavefronts}
  	Assume $\Gamma \cap N^*S=\emptyset$. There is a well-defined sequentially continuous map \[j^*: \mathcal M_\Gamma^{-\infty}(\Gr_{n-k}(V), \Dens(V/E))\to \mathcal M^{-\infty}(\Gr_{m-k}(U), \Dens(U/F)).\]
  	extending the push-forward $J_*$ on measures supported outside $S$. 
  \end{Proposition}
  \begin{Remark}
  	We refer to $j^*$ as the restriction of Crofton distributions.
  \end{Remark}
  
  \proof
  
  Denote $l=n-k$. Consider the partial flag manifold \[F_{l-d,l}=\{(F, E)\in \Gr_{l-d}(V)\times  \Gr_l(V): F\subset E\}\] with the natural projections $\pi_i:F_{l-d,l}\to \Gr_i(V)$ for $i=l,l-d$. As $\pi_{l-d}$ is submersive, $X_U:=\pi_{l-d}^{-1}(\Gr_{l-d}(U))$ is a submanifold of $F_{l-d,l}$. 
  
  Let us check that \[\pi_l^*:\mathcal M_\Gamma^{-\infty}(\Gr_{l}(V), \Dens(V/E))\to \mathcal M^{-\infty}(X_U, \pi_l^*\Dens(V/E))\] is well-defined and sequentially continuous. 
  First, note there is a natural exact sequence \[0\to (E/F)^*\otimes (V/E) \to T_{(F,E)}X_U\to\ F^*\otimes U/F\to 0\] where the inclusion is the differential of the embedding $\Gr_{d}(V/F)\hookrightarrow  X_U$ given by the inverse projection $V\to V/F$, that is $E/F\mapsto (F, E)$; and the surjection is the differential of $\pi_{l-d}$.
  Therefore $\Dens(T_{F,E}X_U)=\Dens^*(E)^k\otimes \Dens(V/E)^d\otimes \Dens(U/F)^{n-(k+d)}$. Moreover, the natural map $f:U/F\rightarrow V/E$ is between spaces of equal dimension, and so induces a map $f^*:\Dens(V/E)\to \Dens(U/F)$. Note that $f$ is an isomorphism if and only if $E\cap U=F$, otherwise $f^*=0$.
  We thus get a natural map $\Dens(T_E\Gr_{l}(V))=\Dens^*(E)^k\otimes\Dens(V/E)^{n-k}\rightarrow \Dens(T_{F,E}X_U)$. This map of line bundles $\pi_l^*|\omega_{\Gr_l(V)}|\to |\omega_{X_U}|$ over $X_U$ combines with the pull-back $\pi_l^*$ to yield the map \[\pi_l^*:\mathcal M^{\infty}(\Gr_{l}(V), \Dens(V/E))\to \mathcal M^{\infty}(X_U, \pi_l^*\Dens(V/E))\]
  denoted by the same letter.
  
  To extend $\pi_l^*$ to the space of distributions as claimed, one has to check that $d\pi_l^*\xi\neq 0$ whenever $(E,\xi)\in \Gamma$, which clearly holds since $\Im(d\pi_l)=TS$. More precisely, if $(F, E)\in X_U$ and $F\subset E\cap U \subset E$ with $\dim E\cap U=l-d+r$ so that $E\in S_r$, then $d\pi_l(T_{E,F}X_U)=T_ES_r$.

  Again noting that for $(F, E)\in X_U$ one has the map $f^*:\Dens(V/E)\to \Dens(U/F)$, we see that the push-forward
  \[(\pi_{l-d})_*: \mathcal M^{-\infty}(X_U, \Dens(V/E))\to \mathcal M^{-\infty}(\Gr_{l-d}(U), \Dens(U/F))\] is well-defined. It remains to define $j^*:=(\pi_{l-d})_*\circ \pi_{l}^*$.
  
  Let us check that $j^*$ extends $J_*$. There is a dense open embedding $i_U:\Gr_l(V)\setminus S\to X_U$ given by $i_U(E)=(E\cap U, E)$. If $\mu\in\mathcal M^{-\infty}(\Gr_l(V), \Dens(V/E))$ and $\supp \mu\cap S=\emptyset$, one  immediately checks that $\pi_l^*\mu=(i_U)_*\mu$, and so  $j^*\mu=(\pi_{l-d})_*\pi_l^*\mu=(\pi_{l-d}\circ i_U)_*\mu=J_*\mu$.
  
  \endproof
  
  \begin{Remark}
  	It follows from the proof that the restriction of a smooth Crofton measure is smooth, as the push-forward of a smooth measure by a submersion is smooth, see e.g. \cite{guillemin_sternberg}.
  \end{Remark}
  
  Let $\widetilde S\subset \AGr_{n-k}(V)$ be the set of translations of elements of $S$, and let $\widetilde J:\AGr_{n-k}(V)\setminus\widetilde S\to \AGr_{m-k}(U)$ be given by $\widetilde J(\widetilde S)= U\cap\widetilde E$. 
  \begin{Lemma}
  	For $f\in C_c^\infty (\AGr_{l-d}(U))$, the function $\widetilde J^*f\in C^\infty (\AGr_l(V)\setminus \widetilde S)$ extends to a compactly supported $\widetilde J^*f\in L^\infty (\AGr_l(V))$.
  \end{Lemma}

  \proof
  Fix a Euclidean structure on $V$. It is obvious that $\widetilde J^*f$ is bounded and smooth outside $\widetilde S$. Considered as a generalized function on $\AGr_l(V)$, it is evidently compactly supported: for a linear subspace $E\in \Gr_l(V)\setminus S$ and $v\perp E$, $(E+v)\cap U=(E\cap U)+ u$ for some $u\in U$. Then there is $e\in E$ such that $u=e+v$, so that $|u|^2=|e|^2+|v|^2$, in particular, $|u|\geq |v|$. Since $f$ is compactly supported, there is $R>0$ such that $f(F+u)=0$ for all $F\in\Gr_{l-d}(U)$, $u\perp F$ and $|u|>R$. Then $\widetilde J^*f(E+v)=0$ for all $|v|>R$.
  \endproof

  \begin{Proposition}\label{prop:Crofton_restriction}
  	Assume $\phi=\Cr(\mu)\in \Val_{k}^{+,KS}(V;\Gamma)$ and $\WF(\mu)\subset\Gamma$. Then 
  	$j^*\phi=\Cr(j^*\mu)$.
  \end{Proposition}
  
  \proof
  We will write $\mu^{tr}$ for the representation of $\mu$ as a translation-invariant distribution on the corresponding affine Grassmannian.
  
  Let us first verify this identity for a smooth valuation $\phi$ given by a smooth Crofton measure $\mu$. Note that $\widetilde J_*\mu^{tr}$ is well-defined in the following sense: given $f\in C^\infty_c(\AGr_{m-k}(U))$, one has $\widetilde J^*f\in L^\infty_c(\AGr_{n-k}(V))$, and so we may set $\langle \widetilde J_*\mu^{tr}, f\rangle=\langle \mu^{tr}, \widetilde J^*f\rangle$. It is also easy to see for smooth $\mu$ that $\widetilde J_*\mu^{tr}=(j^*\mu)^{tr}$. Now
  
  \[ j^*\phi(K)=\phi (jK)=\int_{\AGr_{n-k}(V)} \chi(jK\cap \overline E)d\mu^{tr}(\overline E )= \int_{\AGr_{n-k-d}(U)} \chi (K\cap \overline F)d \widetilde J_*\mu^{tr}(\overline F) \]
  implying $j^*\phi=\Cr(j^* \mu)$. The general result then follows by approximating $\mu_i\to \mu$ by convolving it with an approximate identity on $\GL(V)$, implying $j^*\mu_i\to j^*\mu$ as well as $\Cr(\mu_i)\to \Cr(\mu)$ in $\Val^{KS}(V)$ and thus $j^*\Cr(\mu_i)\to j^*\phi$.
  \endproof
  
  \begin{Remark}
  	One can consider the Klain map of a valuation as the restriction of its Crofton distribution: $\Kl(\Cr(\mu))=i_E^*\mu(E)\in \mathcal M^{-\infty}(\AGr_0(E))^{tr}=\Dens(E)$ where $i_E:E\to V$ is the inclusion.
  \end{Remark}
  
  \subsection{Push-forward by surjection.}\label{sub:fourier}
  
  Again we consider a projection $\pi: V\to W$ with a $d$-dimensional kernel $L$. Define $S_\pi\subset\Gr_{n-k}(V)$ to be the collection of subspaces intersecting $L$ non-transversally.
  
  We could repeat the argument in the previous case. Instead, we will make use of Alesker-Fourier duality, denoted $\mathbb F$, to reduce this case to the previous one. This duality was introduced by Alesker in \cite{alesker_fourier}.
  We will only apply $\mathbb F$ to even valuations, where its description is by far simpler than it is for odd valuations. Namely, $\mathbb F:\Val_k^{+,\infty}(V)\to \Val_{n-k}^{+,\infty}(V^*)\otimes\Dens(V)$ is an isomorphism of Frechet spaces, given in terms of the Klain injection by $\Kl(\mathbb F\phi)(E^\perp)=\Kl(\phi)(E)$. Using the natural extension to spaces of valuations twisted by densities, one has  $\mathbb F^2=\Id$. It extends to an isomorphism of the corresponding spaces of generalized valuations.
  
  For $E\in\Gr_{n-k}(V)$ there is a natural identification $\Dens(V^*/E^\perp)=\Dens^*(E)$. Thus for a dual Crofton distribution $\mu\in \mathcal M^{-\infty}(\Gr_{n-k}(V), \Dens^*(E))$, we may define its Alesker-Fourier transform by \[ {\mathbb  F}\mu:=\perp_*(\mu)\in \mathcal M^{-\infty}(\Gr_{k}(V^*), \Dens(V^*/E)).\] 
  
  \begin{Lemma}\label{lem:Fourier_Crofton_interchange}
  	$\Cr( {\mathbb F}\mu)=\mathbb F\Cr(\mu)$.
  \end{Lemma} 
  \proof
  Using the injectivity of the Klain map, and since $\Kl\circ\Cr$ is the cosine transform while $\Kl\circ\mathbb F=\perp^*\circ \Kl$, this claim reduces to the interchangeability of the cosine transfrom with pull-back by orthogonal complement.
  \endproof
  
  Let $\Gamma\subset T^*\Gr_{n-k}(V)\setminus 0$ be a closed cone. Assume $\Gamma \cap N^*S_\pi=\emptyset$. 
  \begin{Proposition}
  There is then a well-defined sequentially continuous map \[\pi_*: \mathcal M_\Gamma^{-\infty}(\Gr_{n-k}(V), \Dens^*(E))\to \mathcal M^{-\infty}(\Gr_{n-k}(W), \Dens^*(F))\] extending the standard push-forward $\pi_*$ on measures supported outside $S_\pi$.
  \end{Proposition}
  We will simply refer to $\pi_*$ as the projection of Crofton distributions.
  \proof
  
  Consider $\perp:\Gr_{k}(V^*)\to \Gr_{n-k}(V)$. Set $\Gamma'=\perp^*\Gamma\subset T^*\Gr_k(V^*)\setminus 0$, and denote by $S_j:=\perp(S_\pi)$ the $k$-planes intersecting $U=W^\perp$ non-generically. Note that $\Gamma'\cap N^*S_j=\emptyset$, and that  transposed map $\pi^{T}:W^*\to V^*$ is an inclusion. 
  Setting $\pi_*\mu:=\mathbb F\circ (\pi^{T})^*\circ\mathbb F(\mu)$, the statement now follows from Proposition \ref{prop:intersection_wavefronts}. 
  \endproof
  
  We will need the following property of Alesker-Poincar\'{e}  duality proved in \cite{alesker_fourier}
  
  \begin{Proposition}[Alesker]
  	For $\phi\in \Val_k^{\infty}(V)$, it holds that $\mathbb F\pi_*\phi=(\pi ^T)^*\mathbb F\phi$.
  \end{Proposition}
  It follows by approximation that the same holds also for Klain-Schneider continuous valuations, see \cite{bernig_faifman_opq}.
  
  Now we can easily describe the push-forward under projection of a valuation.	
  \begin{Proposition}\label{prop:Crofton_projection}
  	Let $\mu$ be any dual Crofton distribution for $\phi\in \Val_{k}^{+,KS}(V;\Gamma)\otimes \Dens(V^*)$ with $\WF(\mu)\subset\Gamma$. Then $\pi_*\phi=\Cr(\pi_*\mu)$.
  \end{Proposition}
  
  \proof
  Invoking the identity $\mathbb F \circ(\pi ^T)^*\circ\mathbb F(\phi) =\pi_*(\phi)$ for $\phi\in \Val_k^{+,KS}(V)$, this is an immediate consequence of Proposition \ref{prop:Crofton_restriction} and the definition of $\pi_*$.
  
  \endproof

 \end{appendix}

  %\section{Appendix: Generalized valuation as a quasi-completion}

\end{document}